\documentclass[a4paper,12pt]{amsart}
\usepackage{brunnian}
  \usetikzlibrary{arrows}
  \usetikzlibrary{decorations.pathmorphing}
  \usetikzlibrary{external}
\usepackage[labelsep=colon]{caption}
\usepackage{float}
  \restylefloat{figure}
\usepackage{fourier}
\usepackage[margin=3cm]{geometry}
\usepackage{graphicx}
\usepackage{harvard}
\usepackage{hyperref}
\usepackage{mathrsfs}
\usepackage[hang]{subfigure}

\title[On Structure and Organization: An Organizing Principle]{On
  Structure and Organization:\\An Organizing Principle}
\author[N.A. Baas]{Nils A.\ Baas$^\ast$}
\dedicatory{Deparment of Mathematical Sciences, NTNU, N-7491
  Trondheim, Norway\\ \mbox{}\\ (Received 25 June 2012; final version
  received 4 September 2012)}
\thanks{$^\ast$Email: baas@math.ntnu.no}

\DeclareMathOperator{\bond}{Bond}
\DeclareMathOperator{\clu}{Cl}
\DeclareMathOperator{\codim}{codim}
\DeclareMathOperator*{\colim}{colim}
\DeclareMathOperator{\corr}{Corr}
\DeclareMathOperator{\order}{order}
\newcommand{\C}{\mathbb{C}}
\newcommand{\calB}{\mathcal{B}}
\newcommand{\calH}{\mathcal{H}}
\newcommand{\calL}{\mathcal{L}}
\newcommand{\calM}{\mathcal{M}}
\newcommand{\calP}{\mathcal{P}}
\newcommand{\calS}{\mathcal{S}}
\newcommand{\calV}{\mathcal{V}}
\newcommand{\calX}{\mathcal{X}}
\newcommand{\catC}{\mathscr{C}}
\newcommand{\dec}{\mathrm{Dec}}
\newcommand{\R}{\mathbb{R}}
\newcommand{\sets}{\mathrm{Sets}}

\theoremstyle{definition}
\newtheorem{definition}{Definition}
\newtheorem{example}{Example}
\newtheorem*{remark}{Remark}

\begin{document}

\begin{abstract}
  We discuss the nature of structure and organization, and the process
  of making new Things.  Hyperstructures are introduced as binding and
  organizing principles, and we show how they can transfer from one
  situation to another.  A guiding example is the hyperstructure of
  higher order Brunnian rings and similarly structured many-body
  systems.\\

\end{abstract}

\maketitle

\section{Introduction}
\label{sec:intro}
In science and nature we study and utilize collections of objects by
organizing them by relations and patterns in such a way that some
structure emerges.  Objects are bound together to form new objects.
This process may be iterated in order to obtain higher order
collections.  Evolution works along these lines.

When things are being made or constructed it is via binding processes
of some kind.  This seems to be a very general and useful principle
worthy of analyzing more closely.  In other words we are asking for a
general framework in which to study general many-body systems and
their binding patterns as organizing principles.

\section{Examples}
\label{sec:ex}
Let us look at some examples of what we have in mind.

\begin{example}[Links]
  \label{ex:links}
  A link is a disjoint union of embedded circles (or rings) in three
  dimensional space:
  \begin{equation*}
    L \colon \coprod_{i = 1}^n S_i^1 \to \R^3.
  \end{equation*}
  They may be linked in many ways.  Linking is a kind of geometrical
  or topological binding as we see in the following examples.

  \begin{figure}[H]
    \centering
    \subfigure[Hopf link]{
      \begin{tikzpicture}[scale=0.825]
        \draw[knot,double=ring1] (0,0) arc(180:0:1.5);
        \draw[knot,double=ring2] (3,0) circle (1.5);
        \draw[knot,double=ring1] (0,0) arc(-180:0:1.5);
      \end{tikzpicture}
    } \qquad \subfigure[Borromean rings]{
      \begin{tikzpicture}[every path/.style={knot},scale=0.5]
        \pgfmathsetmacro{\csixty}{cos(60)}
        \pgfmathsetmacro{\ssixty}{sin(60)}
        \pgfmathsetmacro{\brscale}{2}
        \draw[double=ring1] (0,0) circle (\brscale);
        \draw[double=ring2] (\brscale,0) circle (\brscale);
        \begin{pgfonlayer}{back}
          \draw[double=ring3] (\brscale *\csixty,\brscale *\ssixty)
            ++(0,\brscale) arc(90:220:\brscale);
          \draw[double=ring3] (\brscale *\csixty,\brscale *\ssixty)
            ++(0,-\brscale) arc(-90:-10:\brscale);
        \end{pgfonlayer}
        \draw[double=ring3] (\brscale *\csixty,\brscale *\ssixty)
          ++(0,\brscale) arc(90:-10:\brscale);
        \draw[double=ring3] (\brscale *\csixty,\brscale *\ssixty)
          ++(0,-\brscale) arc(-90:-140:\brscale);
      \end{tikzpicture}
    } \qquad \subfigure[Brunnian rings]{
      \begin{tikzpicture}[every path/.style={knot},scale=0.5]
        \pgfmathsetmacro{\brscale}{1.8}
        \foreach \brk in {1,2,3} {
          \begin{scope}[rotate=\brk * 120 - 120]
            \colorlet{chain}{ring\brk}
            \brunnianlink{\brscale}{120}
          \end{scope}
        }
      \end{tikzpicture}
    }
    \caption{}
    \label{fig:links}
  \end{figure}

  (Figures in colour are available at:
  \href{http://arxiv.org/pdf/1201.6228.pdf}{http://arxiv.org/pdf/1201.6228.pdf})\\
  
  In \citeasnoun{NS} the linking bonds have been extended to higher order
  links like:

  \begin{figure}[H]
    \centering
    \subfigure[2nd order Brunnian rings]{
      \includegraphics[width=0.4125\linewidth]{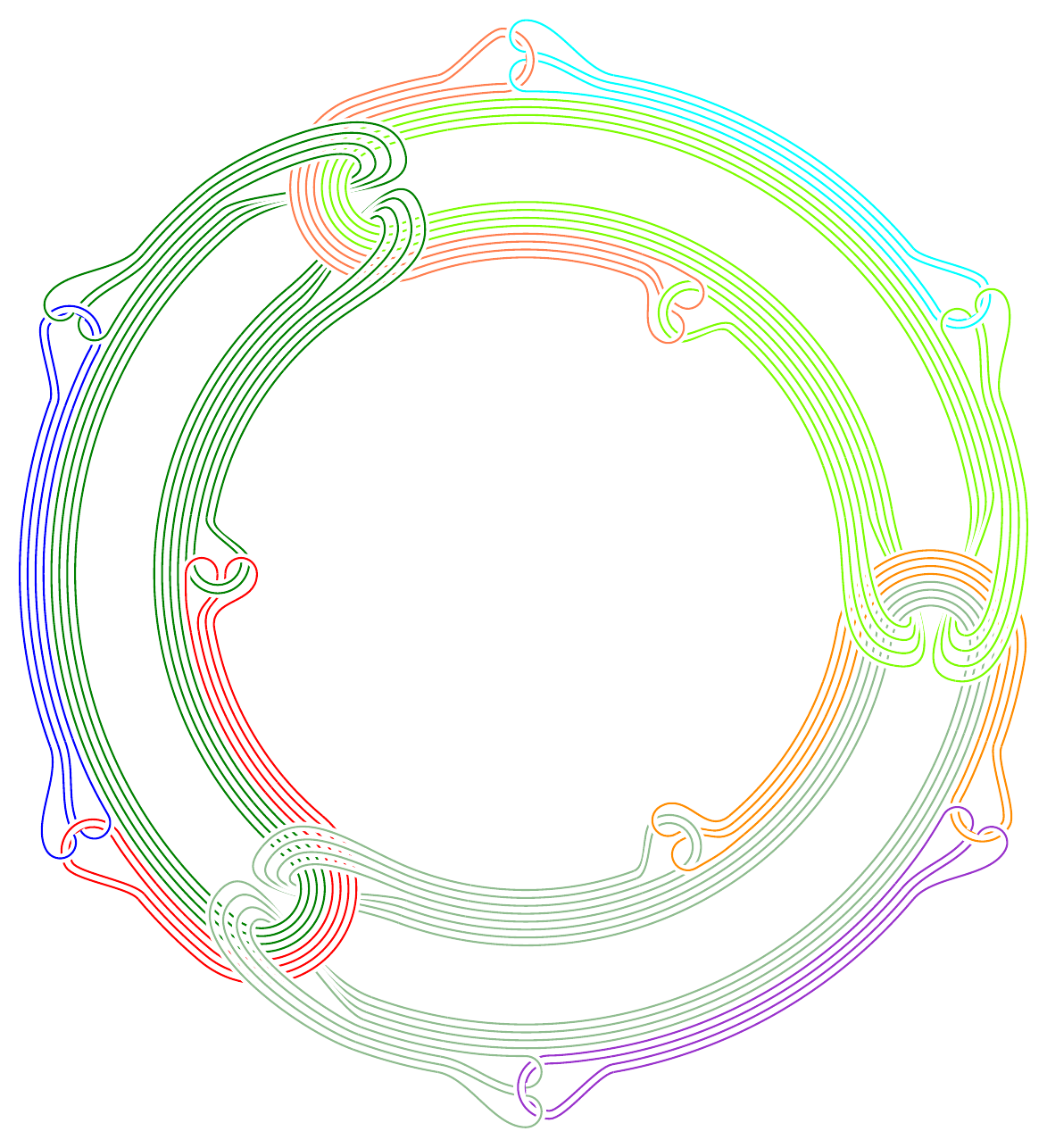}
      \label{fig:2ndBorromean}
    } \qquad \subfigure[3rd order Brunnian rings]{
      \includegraphics[width=0.4125\linewidth]{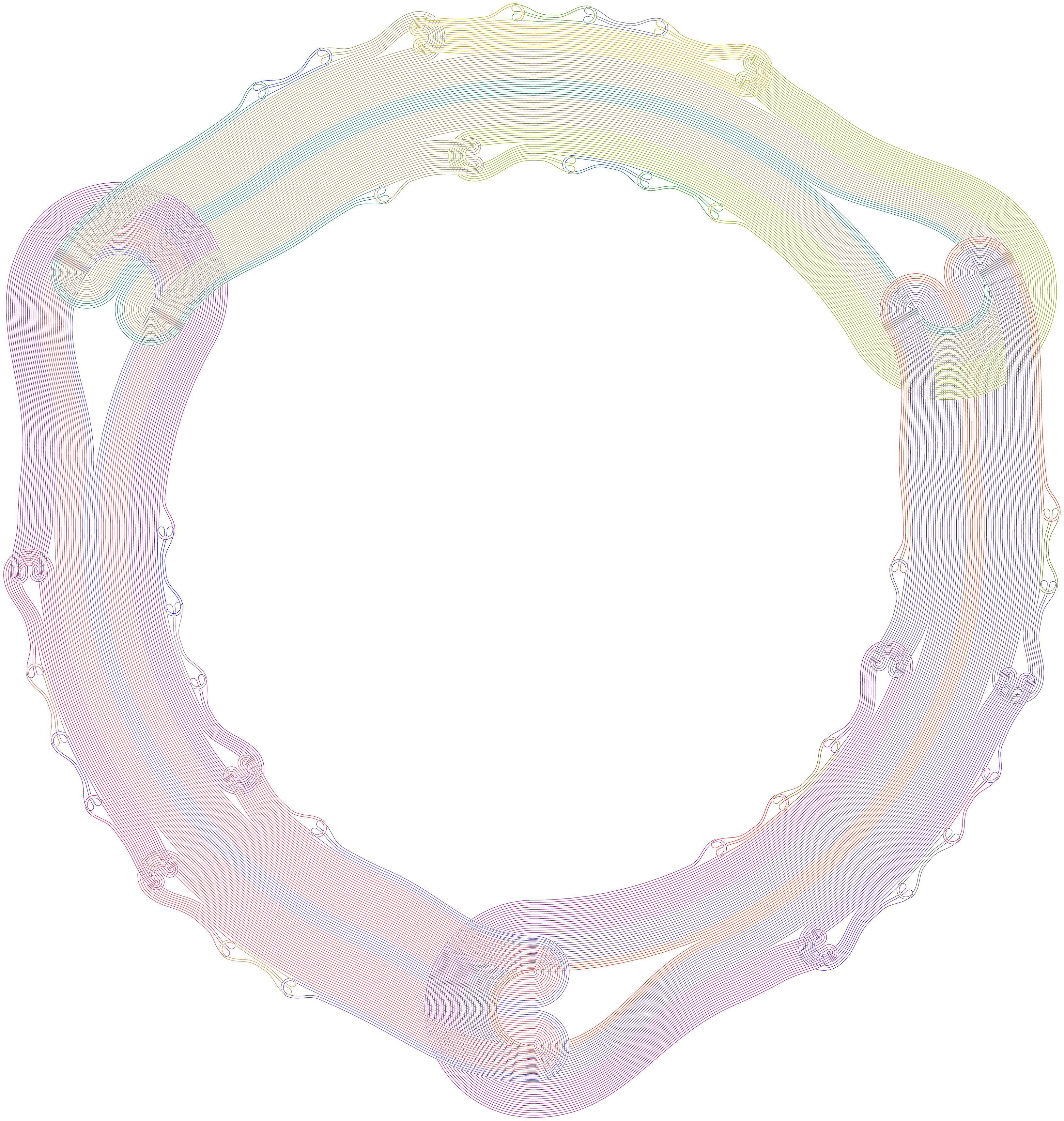}
    }
    \caption{}
    \label{fig:high_links}
  \end{figure}

  In order to iterate this process and study higher order links, it is
  preferable to study embeddings of tori:
  \begin{equation*}
    L \colon \coprod_{i = 1}^n S_i^1 \times D^2 \to S^1 \times D^2
  \end{equation*}

  A second order Brunnian ring binds $9$ circles (rings) together in a
  very subtle way, Figure \ref{fig:2ndBorromean}.  Higher order links
  (links of links of $\ldots$) provide a very good guiding example of
  what a general framework should cover.  For more details, see
  \citeasnoun{NS}.

  In \citeasnoun{NSCS} and \citeasnoun{BS} we discuss possible ways to
  synthesize such binding structures as molecules.
\end{example}

\begin{example}[Many-body states]
  \label{ex:many}
  Efimov (Borromean, Brunnian) states in cold gases are bound states
  of three particles which are not bound two by two.  Hence these
  states are analogous to Borromean and Brunnian rings.  In
  \citeasnoun{NS} we have suggested that this analogy may be extended
  to higher order links and hence suggests higher order versions of
  Efimov states.  For example the second order Brunnian rings
  $2B(3,3)$, see Figure \ref{fig:2ndBorromean}, suggest that there
  should exist bound states of $9$ particles, bound $3$ by $3$ in a
  Brunnian sense, and that these clusters bound together again in a
  higher order sense as in:

  \begin{figure}[H]
    \centering
    \subfigure[$1B(3)$-ring $\sim$ trimers]{
      \begin{tikzpicture}[scale=0.75]
        \draw[very thick] (0,0) circle(2.5cm);
        \draw[very thick, blue] (-0.25,-1.6) circle(0.5cm);
        \draw[very thick, red] (-1.1,1.2) circle(0.5cm);
        \draw[very thick, green] (1.4,0.8) circle(0.5cm);
      \end{tikzpicture}
      \label{fig:trim}
    } \qquad\qquad
    \subfigure[$2B(3,3)$-ring $\sim$ trimers of trimers]{
      \begin{tikzpicture}[scale=0.75]
        \draw[very thick] (0,0) circle(2.5cm);
        
        \begin{scope}[scale=0.3,xshift=-4.25cm,yshift=2.5cm]
          \draw[very thick] (0,0) circle(2.5cm);
          \draw[very thick, blue] (-0.25,-1.6) circle(0.5cm);
          \draw[very thick, red] (-1.1,1.2) circle(0.5cm);
          \draw[very thick, green] (1.4,0.8) circle(0.5cm);
        \end{scope}
        
        \begin{scope}[scale=0.3,xshift=4.5cm,yshift=1.2cm]
          \draw[very thick] (0,0) circle(2.5cm);
          \draw[very thick, blue] (-0.25,-1.6) circle(0.5cm);
          \draw[very thick, red] (-1.1,1.2) circle(0.5cm);
          \draw[very thick, green] (1.4,0.8) circle(0.5cm);
        \end{scope}
        
        \begin{scope}[scale=0.3,xshift=-0.5cm,yshift=-4cm]
          \draw[very thick] (0,0) circle(2.5cm);
          \draw[very thick, blue] (-0.25,-1.6) circle(0.5cm);
          \draw[very thick, red] (-1.1,1.2) circle(0.5cm);
          \draw[very thick, green] (1.4,0.8) circle(0.5cm);
        \end{scope}
      \end{tikzpicture}
      \label{fig:trimtrim}
    }
    \caption{}
    \label{fig:trimtrimtrim}
  \end{figure}

  See also \citeasnoun{NS} for intermediate bound states between a
  trimer and a dimer and a trimer and a two-singleton.  For a
  discussion of higher order Brunnian states, see \citeasnoun{BFJRVZ}.
\end{example}

In general, clustering and higher order clustering of many-body
systems represent a binding mechanism of particle systems ---
parametrizing the particles, in a way.  One may ask for a general
method to describe the binding of particles into higher clusters.

\begin{example}[Clusters and decompositions]
  \label{ex:clust}
  As pointed out in the previous example clusters of objects or data
  represent a binding mechanism between the objects in the cluster.
  The cluster of course depends on various defining criteria.
  Clusters of clusters$\ldots$ represent higher order versions.
  Similarly when we decompose a set or collection we may say that
  elements in the same part of a decomposition are bound together.  In
  this case we get higher bonds as decompositions of
  decompositions$\ldots$
\end{example}

\begin{example}[Mathematical structures and organizations]
  \label{ex:mathstruc}
  We will give a few mathematical examples of how sets are organized
  into structures.

  \begin{itemize}
  \item[a)] Topological spaces.  We organize the ``points'' of a set
    into open sets in such a way that they satisfy the axioms for a
    topology.
  \item[b)] Groups, Algebras, Vector spaces.  The elements are
    organized by certain operations which satisfy the structure
    axioms.
  \item[c)] Manifolds.  Organize the points into open sets and glue
    them together in a prescribed structured way.  Gluing is an
    important example of a geometric binding mechanism.
  \end{itemize}

  In order to form higher order versions of these structures there may
  be many choices, but one way is through higher categories.

  In a higher category of order $n$ one is given objects (e.g.\
  groups, topological spaces,$\ldots$) which are organized by
  morphisms
  \begin{equation*}
    X \xrightarrow{f} Y
  \end{equation*}
  between them.

  Furthermore, there exist morphisms between morphisms ---
  $2$-morphisms:
  \begin{center}
    \begin{tikzpicture}
      \node (X) at (0,0){$X$};
      \node (Y) at (2,0){$Y$};
      
      \path[->,font=\scriptsize]
        (X.north east) edge[bend right=-30]
        node[above]{$f$} (Y.north west)
        (X.south east) edge[bend right=30] node[below]{$g$} (Y.south
        west);
      
      \draw[double equal sign distance,double,-implies] (1,0.3) --
        node[midway,right]{$F$} (1,-0.3);
    \end{tikzpicture}
  \end{center}
  and this continues up to $n$-morphisms between $(n - 1)$-morphisms
  satisfying certain conditions.

  Another type of higher order mathematical binding structure is a
  moduli space.  These are spaces of structures --- for example the
  space of all surfaces of a given genus.
\end{example}

\section{Hyperstructures}
\label{sec:hyper}
We will now introduce some new binding mechanisms of general
collections of objects: physical, chemical, biological, sociological,
abstract and mental.

This organization may bring to light some new and useful structure on
the collections.  We will discuss this in the following, extending the
points of view of \citeasnoun{NS} --- especially in the appendix.

The main concept we will use in order to do this is that of a
\emph{Hyperstructure} as introduced and studied in
\cite{EHH,Hnano,Cog,HAM,NSCS}.

Let us recall the basic construction from \citeasnoun{HAM} and
\citeasnoun{NSCS}.  We start with a set of objects $X_0$ --- our basic
units. To each subset (or families of elements)
\begin{equation*}
  S_0\subset X_0
\end{equation*}
we assign a set of properties or states, $\Omega_0(S_0)$, so
\begin{equation*}
  \Omega_0 \colon \calP(X_0) \to \sets
\end{equation*}
where $\calP(X) = \{A \mid A\subset X\}$ --- the set of subsets ---
the power set, and $\sets$ denotes a suitable set of sets. (In the
language of category theory $\calP(X_0)$ would be considered a
category of subsets, $\sets$ as some category of sets.)  In our
notation here we include properties and states of elements and subsets
of $S_0$ in $\Omega_0 (S_0)$.

Then we want to assign a set of bonds, relations, relationships or
interactions of each subset $S_0$ --- depending on properties and
states.  Here we will just call them bonds.  Let us define
\begin{align*}
  \Gamma_0 &= \{(S_0,\omega_0) \mid S_0\in\calP(X_0),
  \omega_0\in\Omega_0(S_0)\}\\
  B_0 &\colon \Gamma_0\to\sets.
\end{align*}

In our previous notion of hyperstructures the set $X_0$ represents the
systems or agents $(S_i),\Omega_0$ the observables (Obs), $B_0$ the
interactions (Int) and a specific choice of $b_0\in B_0(S_0,\omega_0)$
represents the resultant ``bond'' system giving rise to the next level
of objects --- called $R$ in previous papers, like $S_i$, Obs, Int,
see \cite{EHH,Hnano,Cog,Helvik}.

We will often implicitly assume in the following that given a bond we
know what it binds.  We may require that the set of all bonds of
$(S_0,\omega_0)$ --- $B_0(S_0,\omega_0)$ satisfies the following
condition:\\
\begin{equation*}
  (\ast) \qquad \begin{cases}
    \text{For all } S_0,S'_0,\omega_0,\omega'_0, S_0 \neq S'_0
    \implies B_0(S_0,\omega_0) \cap B_0(S'_0,\omega'_0) = \emptyset\\
    \text{In some cases we may impose the stronger condition} \colon
    B_0 \text{ injective}.
  \end{cases}
\end{equation*}
\mbox{}\\
\noindent (In other situations this is too strong a condition.  For
example if we want to consider colimits as bonds, then the
$\partial_i$ in the following are not well-defined.  The ($\ast$)
condition ensures that the bonds ``know'' what they bind.)\\

\begin{example}[Geometric example]
  \label{ex:geom_ex}
  $S_0 = $ a finite number of manifolds, $\omega_0 = $ property of
  being smooth, put $B_0(S_0,\omega_0) = $ the set of all smooth
  manifolds with boundary equal (isomorphic) to the disjoint union of
  the manifolds in $S_0$.  See Figure \ref{fig:surfaces}.
\end{example}

We now just formalize in a general setting the procedure we described
in \citeasnoun[Section 5]{NSCS}.

Let us form the next level and define:
\begin{equation*}
  X_1 = \{b_0 \mid b_0\in B_0(S_0,\omega_0), S_0\in\calP(X_0),
  \omega_0\in\Omega_0(S_0)\},
\end{equation*}
by definition the image set of $B_0$, and
\begin{center}
  \begin{tikzpicture}
    \node (X1) at (0,1.25){$X_1$};
    \node (calP) at (0,0){$\calP(X_0)$};

    \draw[->] (X1) --
      node[midway,right,font=\scriptsize]{$\partial_0$} (calP);
  \end{tikzpicture}
\end{center}
given by $\partial_0(b_0) = S_0$.

If $B_0$ is injective we have a factorization:
\begin{center}
  \begin{tikzpicture}
    \node (X1) at (0,2.5){$X_1$};
    \node (Gamma) at (1.5,1.25){$\Gamma_0$};
    \node (calP) at (0,0){$\calP(X_0)$};

    \begin{scope}[midway,font=\scriptsize,->]
      \draw (X1) -- node[left]{$\partial_0$} (calP);
      \draw (X1) -- node[above right]{$\partial_0'$} (Gamma);
      \draw (Gamma) -- node[below right]{projection} (calP);
    \end{scope}
  \end{tikzpicture}
\end{center}
Depending on the actual situation we may consider $\partial_0$ and
$\partial_0'$ as boundary maps.

$X_1$ represents the bonds of collections of elements or interactions
in a dynamical context. But the bonds come along with the collection
they bind just as morphisms in mathematics come along with sources and
targets. Similarly at this level we introduce properties and state
spaces and sets of bonds as follows:
\begin{equation*}
  \Omega_1 \colon \calP(X_1) \to \sets
\end{equation*}
($\Omega_1$ then represents the emergent properties as in
\citeasnoun{EHH})

\begin{align*}
  \Gamma_1 &= \{(S_1,\omega_1) \mid S_1\in \calP(X_1),
  \omega_1\in\Omega_1(S_1)\}\\
  B_1 &\colon \Gamma_1\to \sets.
\end{align*}
$B_1$ satisfying the corresponding $(\ast)$ condition.\\

Then we form the next level set:
\begin{equation*}
  X_2 = \{b_1 \mid b_1\in B_1(S_1,\omega_1), S_1\in \calP(X_1),
  \omega_1\in\Omega_1(S_1)\}
\end{equation*} 
and
\begin{center}
  \begin{tikzpicture}
    \node (X2) at (0,1.25){$X_2$};
    \node (calP) at (0,0){$\calP(X_1)$};

    \draw[->] (X2) --
      node[midway,right,font=\scriptsize]{$\partial_1$} (calP);
  \end{tikzpicture}
\end{center}
$\partial_1(b_1)=S_1$. 

We now iterate this procedure up to a general level $N$:

\begin{align*}
  \Omega_{N-1} &\colon\calP(X_{N-1})\to\sets\\
  B_{N-1} &\colon \Gamma_{N-1}\to\sets
\end{align*}
$B_{N - 1}$ satisfying the corresponding $(\ast)$ condition
\begin{equation*}
  X_N = \{b_{N-1} \mid b_{N-1}\in B_{N-1}(S_{N-1},\omega_{N-1}),
  S_{N-1}\in \calP(X_{N-1}), \omega_{N-1}\in\Omega_{N-1}(S_{N-1})\}
\end{equation*}

This is not a recursive procedure since at each level new assignments
take place. The higher order bonds extend the notion of higher
morphisms in higher categories. 

Let us write

\begin{align*}
  \calX &= \{X_0,\ldots,X_N\}\\
  \Omega &= \{\Omega_0,\ldots,\Omega_{N - 1}\}\\
  \calB &= \{B_0,\ldots,B_{N - 1}\}\\
  \partial &= \{\partial_0,\ldots,\partial_{N - 1}\}. 
\end{align*}
where
\begin{center}
  \begin{tikzpicture}
    \node (Xi+1) at (0,1.25){$X_{i + 1}$};
    \node (calP) at (0,0){$\calP(X_i)$};

    \draw[->] (Xi+1) --
      node[midway,right,font=\scriptsize]{$\partial_i$} (calP);
  \end{tikzpicture}
\end{center}
The $\partial_i$'s generalize the source and target maps in the
category theoretical setting, and we think of them as generalized
boundary maps. An Observer mechanism is implicit in the $\Omega_i$'s.
Sometimes one may also want to require maps $I_i \colon X_i \to
X_{i+1}$ or $I_i \colon \calP(X_i) \to X_{i+1}$ --- generalizing the
identity --- such that $\partial_i \circ I_i = \mathrm{id}$.  As for
$\partial_0$ one may also for $\partial_i$ consider $\partial_i'$.

Further mathematical properties to be satisfied will be discussed
elsewhere, for example composition of bonds. We will then also discuss
how to associate a topological space to a hyperstructure --- a
generalized Nerve construction. Bonds may also have internal
structures like topological spaces, manifolds, algebras, vector
spaces, wave functions, fields, etc.

The intuition behind this is:

\begin{equation*} 
  \begin{array}{lcl}
    X_0 & = & \text{objects like atoms, molecules, manifolds, genes,
      organisms}\\
    X_1 & = & \text{bonds} \sim \text{relations, aggregates, clusters,
      interactions, processes,}\ldots\\
    X_2 & = & \text{bonds of bonds} \sim \text{ relations of
      relations, aggregates of aggregates}\\
    & & \text{(possibly overlapping), clusters of clusters,
      interactions of interactions,} \\
    & & \text{processes of processes, }\ldots\text{ etc.} 
  \end{array}
\end{equation*}

\begin{definition}
  \label{def:hyper}
  The system $\calH=(\calX,\Omega,\calB,\partial)$ where the elements
  are related as described, we call \emph{a hyperstructure of order}
  $N$.
\end{definition}

Sometimes one may want to organize a set of agents for a specific
purpose.  One way to do this is to put a hyperstructure on it
organizing the agents to fulfill a given goal.  This applies to both
concrete physical and abstract situations.

An example of this is the procedure whereby we organize molecules (or
abstract topological bonds) into rings, $2$-rings,$\ldots$, $n$-rings
representing new topological structures \citeasnoun{BS}.

In many cases it is natural to view the bonds as geometric or
topological spaces. For example if a surface $F$ has three circles as
boundary components, see Figure \ref{fig:surfaces}, $\partial F = S_1
\cup S_2 \cup S_3$, we may say that $F$ is a geometric bond of the
circles. Clearly there may be many. This is in analogy with chemical
bonds.

Furthermore, if we have a manifold $B$ such that its boundary is
\begin{equation*}
  \partial B = A_1 \cup \cdots \cup A_n \cup \tilde{A},
\end{equation*}
see Figure \ref{fig:boundary}, we may say that $B$ is a bond of
$A_1,\ldots,A_n$. Even more general is the following situation: given
a topological space of some kind $Y$ and let $Z_1,\ldots,Z_n$ be
subspaces of $Y$. Then we say that $Y$ is a bond of $Z_1,\ldots,Z_n$,
see Figure \ref{fig:slices}, thinking of $Z_1,\ldots,Z_n$ as ``the
boundary'' of $Y$.

In a hyperstructure of higher order we may let $Y$ represent a top
level bond, then the $Z_i$'s will represent bonds of other spaces:
\begin{align*}
  Y &= B_0(Z_1,\ldots,Z_n)\\
  Z_i &= B_{1,i}(W_{1,i},\ldots,W_{n,i})
\end{align*}
etc., see Figure \ref{fig:towers}.

This slightly extends the pictures of \citeasnoun{HAM} and represent
what we could call a \emph{geometric hyperstructure}.  This concept
relates to topological quantum field theory and will be studied in a
separate paper, see Section \ref{sec:multi} for some further remarks.

Hyperstructures and higher order bonds may be viewed as a huge
extension of cobordisms (manifolds with boundary) and chemicals bonds
and interlockings.

Furthermore, the whole scheme of thinking may be applied to
interactions and ways of connecting and interlocking people, groups of
people, social and economic systems and organizations, organisms,
genes, data, etc.  For some of these aspects, see \citeasnoun{HAM}.

By describing all such types of systems by means of hyperstructures
one may create entirely new structures which may be useful both in old
and new contexts.

\begin{figure}[H]
  \centering
  \begin{tikzpicture}[scale=0.5]
    \foreach \x in {2,6,10}{
      \draw (\x,0) arc(0:-180:1cm and 0.5cm);
      \draw[dashed] (\x,0) arc(0:180:1cm and 0.5cm);
    }

    \draw (0,0) .. controls (0,6) and (10,6) ..  (10,0);
    \draw (2,0) .. controls (2,2) and (4,2) .. (4,0);
    \draw (6,0) .. controls (6,2) and (8,2) .. (8,0);

    \node at (1,-1.125){$S_1$};
    \node at (5,-1.125){$S_2$};
    \node at (9,-1.125){$S_3$};
    \node at (5,3){$F$};

    \node at (12,2){or};

    \begin{scope}[xshift=14cm]
      \foreach \x in {2,6,10}{
        \draw (\x,0) arc(0:-180:1cm and 0.5cm);
        \draw[dashed] (\x,0) arc(0:180:1cm and 0.5cm);
      }
      
      \draw (0,0) .. controls (0,6) and (10,6) ..  (10,0);
      \draw (2,0) .. controls (2,2) and (4,2) .. (4,0);
      \draw (6,0) .. controls (6,2) and (8,2) .. (8,0);
      
      \node at (1,-1.125){$S_1$};
      \node at (5,-1.125){$S_2$};
      \node at (9,-1.125){$S_3$};
      \node at (3,3){$F'$};
    \end{scope}

    \begin{scope}[xshift=17.75cm,yshift=2.35cm]
      \draw (0,0) .. controls (0,-0.5) and (1.5,-0.5) .. (1.5,0);
      \draw (0.125,-0.195) .. controls (0.125,0.25) and (1.375,0.25)
        .. (1.375,-0.195);
    \end{scope}

    \begin{scope}[xshift=20cm,yshift=2.75cm]
      \draw (0,0) .. controls (0,-0.5) and (1.5,-0.5) .. (1.5,0);
      \draw (0.125,-0.195) .. controls (0.125,0.25) and (1.375,0.25)
        .. (1.375,-0.195);
    \end{scope}
  \end{tikzpicture}
  \caption{}
  \label{fig:surfaces}
\end{figure}

\begin{figure}[H]
  \centering
  \begin{tikzpicture}[scale=0.5]
    \foreach \x in {2,6,14,18}{
      \draw (\x,0) arc(0:-180:1cm and 0.5cm);
      \draw[dashed] (\x,0) arc(0:180:1cm and 0.5cm);
    }

    \draw[dashed] (8,0) -- (10,0);

    \draw (0,0) .. controls (0,8) and (18,8) ..  (18,0);
    \foreach \a in {2,6,10,14}{
      \draw (\a,0) .. controls(\a,2) and (\a+2,2) .. (\a+2,0);
    }

    \node at (1,-1.125){$A_1$};
    \node at (5,-1.125){$A_2$};
    \node at (13,-1.125){$A_n$};
    \node at (17,-1.125){$\tilde{A}$};
    \node at (9,3){$B$};
  \end{tikzpicture}
  \caption{}
  \label{fig:boundary}
\end{figure}

\begin{figure}[H]
  \centering
  \begin{tikzpicture}[scale=0.7]
    \draw[rotate=15] (0.24396864,-1.4967076) .. controls
      (0.49101374,-2.4657116) and (0.665201,-1.8994577)
      .. (1.5439687,-2.3767076) .. controls (2.4227362,-2.8539574) and
      (3.017664,-3.2517743) .. (4.0039687,-3.4167075) .. controls
      (4.9902735,-3.581641) and (5.629609,-3.3427677)
      .. (6.6239686,-3.2367077) .. controls (7.6183286,-3.1306474) and
      (8.599862,-3.4573336) .. (9.383968,-2.8367076) .. controls
      (10.168076,-2.2160816) and (9.585438,-2.3502257)
      .. (9.943969,-1.4167076) .. controls (10.302499,-0.48318952) and
      (10.588963,-0.26462105) .. (10.743969,0.7232924) .. controls
      (10.898975,1.7112058) and (11.202379,1.6006097)
      .. (10.663969,2.4432924) .. controls (10.125558,3.2859752) and
      (9.734352,3.3049438) .. (8.743969,3.4432924) .. controls
      (7.753585,3.581641) and (5.849801,3.2754364)
      .. (5.0439687,2.6832924) .. controls (4.2381363,2.0911484) and
      (5.064207,1.7455099) .. (4.2239685,1.2032924) .. controls
      (3.3837304,0.661075) and (2.3561227,1.5378368)
      .. (1.5239687,0.9832924) .. controls (0.69181454,0.428748) and
      (-0.003076456,-0.5277036) .. (0.24396864,-1.4967076);

    \begin{scope}[xshift=2cm,yshift=-1.5cm]
      \draw (0,0) .. controls (1,0.5) and (-0.5,1.25) .. (0,2);
      \node[right] at (0.25,0){$Z_1$};
    \end{scope}

    \begin{scope}[xshift=4cm,yshift=-1cm]
      \draw (0,0) .. controls (1,0.5) and (-0.5,1.25) .. (0,2);
      \node[right] at (0.25,0){$Z_2$};
    \end{scope}

    \begin{scope}[xshift=6cm,yshift=-0.5cm,rotate=20]
      \draw[dotted] (0,1) -- (1,1);
    \end{scope}

    \begin{scope}[xshift=8cm,yshift=0cm]
      \draw (0,0) .. controls (1,0.5) and (-0.5,1.25) .. (0,2);
      \node[right] at (0.25,0){$Z_n$};
    \end{scope}

    \node at (7.25,3.5){$Y$};
  \end{tikzpicture}
  \caption{}
  \label{fig:slices}
\end{figure}

\begin{figure}[H]
  \centering
  \begin{tikzpicture}[scale=0.7]
    \draw[rotate=15] (0.24396864,-1.4967076) .. controls
      (0.49101374,-2.4657116) and (0.665201,-1.8994577)
      .. (1.5439687,-2.3767076) .. controls (2.4227362,-2.8539574) and
      (3.017664,-3.2517743) .. (4.0039687,-3.4167075) .. controls
      (4.9902735,-3.581641) and (5.629609,-3.3427677)
      .. (6.6239686,-3.2367077) .. controls (7.6183286,-3.1306474) and
      (8.599862,-3.4573336) .. (9.383968,-2.8367076) .. controls
      (10.168076,-2.2160816) and (9.585438,-2.3502257)
      .. (9.943969,-1.4167076) .. controls (10.302499,-0.48318952) and
      (10.588963,-0.26462105) .. (10.743969,0.7232924) .. controls
      (10.898975,1.7112058) and (11.202379,1.6006097)
      .. (10.663969,2.4432924) .. controls (10.125558,3.2859752) and
      (9.734352,3.3049438) .. (8.743969,3.4432924) .. controls
      (7.753585,3.581641) and (5.849801,3.2754364)
      .. (5.0439687,2.6832924) .. controls (4.2381363,2.0911484) and
      (5.064207,1.7455099) .. (4.2239685,1.2032924) .. controls
      (3.3837304,0.661075) and (2.3561227,1.5378368)
      .. (1.5239687,0.9832924) .. controls (0.69181454,0.428748) and
      (-0.003076456,-0.5277036) .. (0.24396864,-1.4967076);

    \begin{scope}[xshift=2cm,yshift=-1.5cm,scale=0.75,rotate=10]
      \draw (0,0) rectangle(1,3);
      \draw (0,1) rectangle(1,2);
      
      \node[right] at (1.125,0){$Z_1$};
    \end{scope}

    \begin{scope}[xshift=4cm,yshift=-1cm,scale=0.75,rotate=10]
      \draw (0,0) rectangle(1,3);
      \draw (0,1) rectangle(1,2);

      \node[right] at (1.125,0){$Z_2$};
    \end{scope}

    \begin{scope}[xshift=6cm,yshift=-0.5cm,rotate=10]
      \draw[dotted] (0,1) -- (1,1);
    \end{scope}

    \begin{scope}[xshift=8cm,yshift=0cm,scale=0.75,rotate=10]
      \draw (0,0) rectangle(1,3);
      \draw (0,1) rectangle(1,2);

      \node[right] at (1.125,0){$Z_n$};
    \end{scope}

    \node at (7,3.5){$Y$};
  \end{tikzpicture}
  \caption{}
  \label{fig:towers}
\end{figure}

We could also call a hyperstructure a binding structure since it
really binds the elements of a collection.  To make the notion simpler
we will suppress the states in the following, and we will express the
associated binding mechanisms of collections as simple as possible.
At the end of \citeasnoun{HAM} we offer a categorical version which is
more restrictive.

Further examples of hyperstructures:
\begin{itemize}
\item[a)] \emph{Higher order links} as in Example \ref{ex:links}.
  Here the starting set $X_0$ is a collection of rings, the observed
  state is circularity and the bonds are the links.  Then one
  observes ``circularity'' (or embedding in a torus) of the Brunnian
  links and continues the process by forming rings of rings $\ldots$
  as described in \citeasnoun{NS}.\\
\item[b)] \emph{Hierarchical clustering} and multilevel decompositions
  are typical examples of hyperstructures.  The bonds are given by
  level-wise membership as in Example \ref{ex:many}.  We should point
  out that hyperstructures encompass and are far more sophisticated
  and subtle than hierarchies.\\
\item[c)] \emph{Higher categories} are examples of hyperstructures as
  well where the higher morphisms are interpreted as higher bonds.\\
\item[d)] \emph{Compositions}.  Often collections of objects and data
  may be organized into a composition of mappings
  \begin{equation*}
    S_1 \xleftarrow{\varphi_1} S_2 \xleftarrow{\varphi_2}  \cdots
    \xleftarrow{\varphi_{n - 1}} S_n
  \end{equation*}
  where we may think of $s_i\in S_i$ as a bond of the elements in
  $\varphi_i^{-1}(s_i)$, and $s_{i + 1}\in \varphi^{-1}(s_i)$ is again
  a bond of the elements in $\varphi_{i + 1}^{-1}(s_{i + 1})$, etc.
  See \citeasnoun{HAM} and references therein.  Similarly one may say
  that a subset (or space if we have more structure) $Z_i \subset S_i$
  is a bond of subsets in $\varphi_i^{-1}(Z_i)$.  Composition models
  of hyperstructures are particularly interesting when the $S_i{}'s$
  and mappings have more structure, for example being smooth manifolds
  and smooth mappings.  In that case there is an interesting stability
  theory, see \citeasnoun{HAM} and references therein.\\
\item[e)] \emph{Higher Order Cobordisms}. In geometry and topology we
  consider kinds of generalized surfaces in arbitrary dimensions
  called manifolds.  These may be smooth and have various additional
  structures.  Amongst manifolds there is a very important notion of
  cobordism, and we will illustrate how cobordisms of manifolds with
  boundaries and corners are important as bonds.

  Two manifolds $A_1$ and $A_2$ (with or without boundary) of
  dimension $n$ are cobordant if and only if there exists an $(n + 1)$
  dimensional manifold $B$ such that
  \begin{equation*}
    \partial B = (A_1 \cup A_2) \cup \hat{A}
  \end{equation*}
  $\partial$ stands for boundary and $\cup$ means glued together along
  the common boundary
  \begin{equation*}
    \partial(A_1 \cup A_2) = \partial \hat{A},
  \end{equation*}
  see Figure \ref{fig:cob1}. 

  In this paper we are interested in studies of structures etc.  So
  let us see what cobordisms of cobordisms or more generally bonds of
  bonds mean in this geometric setting.

  Let $B_1$ be a bond between $A_1$ and $A_2$ and let $B_2$ be a bond
  between $\tilde{A}_1$ and $\tilde{A}_2$.  Then $C$ is bond
  (cobordism) between $B_1$ and $B_2$ if and only if
  \begin{equation*}
    \partial C = (B_1 \cup B_2) \cup \hat{B}
  \end{equation*}
  where $\cup$ means glued together along common boundary:
  $\partial(B_1 \cup B_2) = \partial \hat{B}, C$ is of dimension $(n +
  2)$ and $\hat{B}$ of dimension $(n+1)$, see Figure \ref{fig:cob2}.

  Furthermore, a third order bond between $C_1$ and $C_2$ will be
  given by an $(n+3)$ manifold $D$ such that
  \begin{equation*}
    \partial D = (C_1 \cup C_2) \cup \hat{C}
  \end{equation*}
  etc., see Figure \ref{fig:cube}. 

  In this formal description we have just considered two
  ``components'' in the boundary, hence a cobordism is then considered
  as a bond between two parts like a morphism in a category.  But
  clearly the geometry extends to any finite number of components,
  hence we consider a cobordism as the prototype of a geometric bond
  between several objects:

  \begin{equation*}
    B(A_1, \ldots ,A_n) \quad \text{if and only if} \quad \partial B =
    (A_1 \cup \cdots \cup A_n) \cup \hat{A}.
  \end{equation*}

  Mathematically this requires that we study manifolds with decomposed
  boundaries, whose boundary components again are decomposed, etc. (as
  introduced and studied in \citeasnoun{Sing}) or manifolds with
  higher order corners (corners of corners etc.), see Figure
  \ref{fig:cob3}.

  Hyperstructures seem like the correct mathematical structure to
  describe this situation. 

  \begin{figure}[H]
    \centering
    \begin{tikzpicture}
      \begin{scope}[scale=0.7]
        \draw[fill=black] (0,2) circle(0.075cm);
        \draw[fill=black] (0,0) circle(0.075cm);
        
        \node[yshift=0.5cm] at (0,2){$A_2$};
        \node[yshift=-0.5cm] at (0,0){$A_1$};
        
        \draw (0,2) -- node[right]{$B$} (0,0);
        
        \begin{scope}[xshift=2cm]
          \draw (0,2) ellipse(0.5cm and 0.25cm);
          \draw (0.5,0) arc(0:-180:0.5cm and 0.25cm);
          \draw[dashed] (0.5,0) arc(0:180:0.5cm and 0.25cm);
          
          \draw (-0.5,0) -- (-0.5,2);
          \draw (0.5,0) -- (0.5,2);
          
          \node[yshift=0.6cm] at (0,2){$A_2$};
          \node[yshift=-0.6cm] at (0,0){$A_1$};
          \node at (0,1){$B$};
        \end{scope}
        
        \begin{scope}[xshift=3.5cm]
          \draw (0.5,2) ellipse (0.5cm and 0.25cm);
          \draw (2.5,2) ellipse (0.5cm and 0.25cm);
          \draw (4.5,2) ellipse (0.5cm and 0.25cm);
          \draw (2,0) arc(0:-180:0.5cm and 0.25cm);
          \draw[dashed] (2,0) arc(0:180:0.5cm and 0.25cm);
          \draw (4,0) arc(0:-180:0.5cm and 0.25cm);
          \draw[dashed] (4,0) arc(0:180:0.5cm and 0.25cm);
          
          \draw (1,2) .. controls (1,1.5) and (2,1.5) .. (2,2);
          \draw (3,2) .. controls (3,1.5) and (4,1.5) .. (4,2);
          
          \draw (2,0) .. controls (2,0.5) and (3,0.5) .. (3,0);
          
          \draw (0,2) .. controls (0,1) and (1,1) .. (1,0);
          \draw (5,2) .. controls (5,1) and (4,1) .. (4,0);
          
          \node[yshift=0.6cm] at (2.5,2){$A_2$};
          \node[yshift=-0.6cm] at (2.5,0){$A_1$};
          \node at (3,1){$B$};
          
          \begin{scope}[xshift=1.25cm,yshift=1cm,scale=0.5]
            \draw (0,0) .. controls (0,-0.5) and (1.5,-0.5)
              .. (1.5,0);
            \draw (0.125,-0.195) .. controls (0.125,0.25) and
              (1.375,0.25) .. (1.375,-0.195);
          \end{scope}
        \end{scope}

        \begin{scope}[xshift=9.5cm]
          \draw (0.5,2) ellipse (0.5cm and 0.25cm);
          \draw (2.5,2) ellipse (0.5cm and 0.25cm);
          \draw (4.5,2) ellipse (0.5cm and 0.25cm);
          \draw (6.5,2) ellipse (0.5cm and 0.25cm);
          \draw (8.5,2) ellipse (0.5cm and 0.25cm);
          
          \foreach \a in {1,3,5,7}{
            \draw (\a,2) .. controls (\a,1.5) and (\a+1,1.5)
              .. (\a+1,2);
          }
          
          \draw (0,2) .. controls (0,-1) and (9,-1) .. (9,2);

          \node[yshift=0.6cm] at (4.5,2){$A_1 \cup A_2$};
          \node at (5.5,0.8){$B$};

          \begin{scope}[xshift=3cm,yshift=0.8cm,scale=0.5]
            \draw (0,0) .. controls (0,-0.5) and (1.5,-0.5)
              .. (1.5,0);
            \draw (0.125,-0.195) .. controls (0.125,0.25) and
              (1.375,0.25) .. (1.375,-0.195);
          \end{scope}
        \end{scope}
      \end{scope}

      \begin{scope}[scale=0.7,xshift=2.25cm,yshift=-6.5cm]
        \begin{scope}
          \draw (0,1) arc(90:270:0.25cm and 0.5cm);
          \draw[dashed] (0,1) arc(90:-90:0.25cm and 0.5cm);
          \draw (0,3) arc(90:270:0.25cm and 0.5cm);
          \draw[dashed] (0,3) arc(90:-90:0.25cm and 0.5cm);
          \draw (5,1.5) ellipse(0.25cm and 0.5cm);
          \draw (3,3.5) ellipse(0.5cm and 0.25cm);
          
          \draw (0,0) .. controls (2.5,0) and (2.5,1) .. (5,1);
          \draw (0,1) .. controls (0.5,1) and (0.5,2) .. (0,2);
          \draw (0,3) .. controls (1,3) and (2.25,2) .. (2.5,3.5);
          \draw (3.5,3.5) .. controls (3.75,2) .. (5,2);
          
          \begin{scope}[xshift=1.5cm,yshift=1.5cm,scale=0.5]
            \draw (0,0) .. controls (0,-0.5) and (1.5,-0.5)
              .. (1.5,0);
            \draw (0.125,-0.195) .. controls (0.125,0.25) and
              (1.375,0.25) .. (1.375,-0.195);
          \end{scope}

          \node at (-1,1.5){$A_1$};
          \node at (6,1.5){$A_2$};
          \node at (3,4.25){$\hat{A}$};
          \node at (3,2.25){$B$};
        \end{scope}

        \begin{scope}[xshift=10cm,yshift=0.25cm]
          \draw (0,0) rectangle(3,3);
          \node[left] at (0,1.5){$A_1$};
          \node[right] at (3,1.5){$A_2$};
          \node at (1.5,1.5){$B$};
          \node[above] at (1.5,3){$\hat{A}$};
          \node[below] at (1.5,0){$\hat{A}$};
        \end{scope}
      \end{scope}
    \end{tikzpicture}
    \caption{}
    \label{fig:cob1}
  \end{figure}

  \begin{figure}[H]
    \centering
    \begin{tikzpicture}
      \begin{scope}[scale=0.575]
        \begin{scope}
          \foreach \y in {-2.5,0.5,2.5,4.5}{
            \draw (5,\y) ellipse(0.25cm and 0.5cm);
          }

          \draw (5,1) .. controls (4,1) and (4,2) .. (5,2);
          \draw (5,3) .. controls (4,3) and (4,4) .. (5,4);

          \draw (5,0) .. controls (3.5,0) and (3.5,-2) .. (5,-2);
          \draw (5,5) .. controls (0,5) and (0,-3) .. (5,-3);

          \node[left] at (1,1){$\hat{B}$};

          \begin{scope}[xshift=2.5cm,yshift=1.5cm,scale=0.5]
            \draw (0,0) .. controls (0,-0.5) and (1.5,-0.5)
              .. (1.5,0);
            \draw (0.125,-0.195) .. controls (0.125,0.25) and
              (1.375,0.25) .. (1.375,-0.195);
          \end{scope}
        \end{scope}

        \begin{scope}[xshift=6cm]
          \foreach \y in {1,3}{
            \draw (0,\y) arc(90:270:0.25cm and 0.5cm);
            \draw[dashed] (0,\y) arc(90:-90:0.25cm and 0.5cm);
            \draw (0,\y) .. controls (1,\y) and (1,\y+1) .. (0,\y+1);
          }
          \draw (0,5) arc(90:270:0.25cm and 0.5cm);
          \draw[dashed] (0,5) arc(90:-90:0.25cm and 0.5cm);

          \draw (5,1.5) ellipse(0.25cm and 0.5cm);
          \draw (5,3.5) ellipse(0.25cm and 0.5cm);
          \draw (5,2) .. controls (4,2) and (4,3) .. (5,3);

          \draw (0,0) .. controls (2.5,0) and (2.5,1) .. (5,1);
          \draw (0,5) .. controls (2.5,5) and (2.5,4) .. (5,4);

          \node[above] at (0,5){$A_1$};
          \node[above] at (2.5,5){$B_1$};
          \node[above] at (5,5){$A_2$};

          \begin{scope}[xshift=2cm,yshift=2.5cm,scale=0.5]
            \draw (0,0) .. controls (0,-0.5) and (1.5,-0.5)
              .. (1.5,0);
            \draw (0.125,-0.195) .. controls (0.125,0.25) and
              (1.375,0.25) .. (1.375,-0.195);
          \end{scope}

          \node at (2.5,-1){$C$};

          \draw (0,-2) arc(90:270:0.25cm and 0.5cm);
          \draw[dashed] (0,-2) arc(90:-90:0.25cm and 0.5cm);
          \draw (5,-2.5) ellipse(0.25cm and 0.5cm);
          \draw (0,-2) -- (5,-2);
          \draw (0,-3) -- (5,-3);

          \node[below] at (0,-3){$\hat{A}_1$};
          \node[below] at (2.5,-3.15){$B_2$};
          \node[below] at (5,-3){$\hat{A}_2$};

          \begin{scope}[xshift=2cm,yshift=-2.5cm,scale=0.5]
            \draw (0,0) .. controls (0,-0.5) and (1.5,-0.5)
              .. (1.5,0);
            \draw (0.125,-0.195) .. controls (0.125,0.25) and
              (1.375,0.25) .. (1.375,-0.195);
          \end{scope}
        \end{scope}

        \begin{scope}[xshift=12cm]
          \foreach \y in {-2,2,4}{
            \draw (0,\y) arc(90:270:0.25cm and 0.5cm);
            \draw[dashed] (0,\y) arc(90:-90:0.25cm and 0.5cm);
          }

          \draw (0,3) .. controls (1,3) and (1,2) .. (0,2);
          \draw (0,1) .. controls (1.5,1) and (1.5,-2) .. (0,-2);
          \draw (0,4) .. controls (5,4) and (5,-3) .. (0,-3);

          \node[right] at (4,1){$\hat{B}$};

          \begin{scope}[xshift=1.5cm,yshift=0.75cm,scale=0.5]
            \draw (0,0) .. controls (0,-0.5) and (1.5,-0.5)
              .. (1.5,0);
            \draw (0.125,-0.195) .. controls (0.125,0.25) and
              (1.375,0.25) .. (1.375,-0.195);
          \end{scope}

          \begin{scope}[xshift=2cm,yshift=-0.25cm,scale=0.5]
            \draw (0,0) .. controls (0,-0.5) and (1.5,-0.5)
              .. (1.5,0);
            \draw (0.125,-0.195) .. controls (0.125,0.25) and
              (1.375,0.25) .. (1.375,-0.195);
          \end{scope}
        \end{scope}
      \end{scope}

      \begin{scope}[xshift=11cm,yshift=-0.5cm,scale=0.7]
        \draw (0,0) rectangle (3,3);

        \node[below left] at (0,0){$\hat{A}_1$};
        \node[left] at (0,1.5){$\hat{B}$};
        \node[above left] at (0,3){$A_1$};
        \node[below] at (1.5,0){$B_2$};
        \node at (1.5,1.5){$C$};
        \node[above] at (1.5,3){$B_1$};
        \node[below right] at (3,0){$\hat{A}_1$};
        \node[right] at (3,1.5){$\hat{B}$};
        \node[above right] at (3,3){$A_1$};

        \foreach \x in {0,3}{
          \foreach \y in {0,3}{
            \draw[fill=black] (\x,\y) circle(0.075cm);
          }
        }
      \end{scope}
    \end{tikzpicture}
    \caption{}
    \label{fig:cob2}
  \end{figure}
  
  \begin{figure}[H]
    \centering
    \begin{tikzpicture}
      \draw (1,0) rectangle(3,2);
      \draw (1,0) -- (0,0.5) -- (0,2.5) -- (2,2.5) -- (3,2);
      \draw (1,2) -- (0,2.5);
      \draw[dashed] (0,0.5) -- (2,0.5) -- (3,0);
      \draw[dashed] (2,2.5) -- (2,0.5);

      \draw[->] (-0.5,2) node[above left]{$C_1$} -- (0.5,1);
      \draw[->] (4,1.75) node[yshift=0.05cm,right]{$C_2$} --
        (2.5,1.25);
      \draw[->] (3,2.75) node[yshift=0.1cm,right]{$D$} --
        (1.75,2.15);
    \end{tikzpicture}
    \caption{$D$ is schematically represented by a cube}
    \label{fig:cube}
  \end{figure}
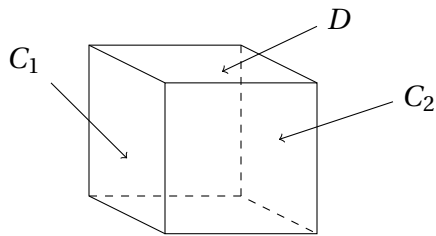
  
  \begin{figure}[H]
    \centering
    \begin{tikzpicture}
      \draw (0.5,2) ellipse (0.5cm and 0.25cm);
      \draw (2.5,2) ellipse (0.5cm and 0.25cm);
      \draw (4.5,2) ellipse (0.5cm and 0.25cm);
      \draw[dashed] (6,2) -- (7,2);
      \draw (8.5,2) ellipse (0.5cm and 0.25cm);
      
      \foreach \a in {1,3,5,7}{
        \draw (\a,2) .. controls (\a,1.5) and (\a+1,1.5) .. (\a+1,2);
      }
          
      \draw (0,2) .. controls (0,-1) and (9,-1) .. (9,2);
      
      \node[yshift=0.6cm] at (0.5,2){$A_1$};
      \node[yshift=0.6cm] at (2.5,2){$A_2$};
      \node[yshift=0.6cm] at (4.5,2){$A_3$};
      \node[yshift=0.6cm] at (8.5,2){$A_n$};
      \node at (5.5,0.8){$B$};
      \node at (10.5,1.5){$, \quad \hat{A} = \emptyset$};
      
      \begin{scope}[xshift=3cm,yshift=0.8cm,scale=0.5]
        \draw (0,0) .. controls (0,-0.5) and (1.5,-0.5) .. (1.5,0);
        \draw (0.125,-0.195) .. controls (0.125,0.25) and (1.375,0.25)
          .. (1.375,-0.195);
      \end{scope}
    \end{tikzpicture}
    \caption{}
    \label{fig:cob3}
  \end{figure}
  
  In the framework we have introduced the geometric examples in the
  figures correspond to:
  \begin{equation*}
    X_0 = \{\text{the set of circles in some high dimensional space}\}
  \end{equation*}
  no states, $\Omega_0=\emptyset$. 

  $S_0\in \calP(X_0)$ means that $S_0$ is a disjoint union of
  circles. $b_0\in B_0(S_0)$ is then given by a surface having the
  circles of $S_0$ as its boundary
  \begin{equation*}
    X_1 = \{\text{the set of surfaces with boundary equal the union of
      circles}\}
  \end{equation*}
  $S_1\in \calP(X_1)$ is a disjoint union of such surfaces.

  $B_1(S_1)$ is then given by a $3$ dimensional manifold having the
  surfaces of $S_1$ as parts of its boundary, but possibly glued
  together along common boundaries with additional parts --- the
  $\hat{B}$'s.  For more details on hyperstructured glueing and
  decomposition processes, see \citeasnoun{Sing} and \citeasnoun{BCR}.

  In this way it goes on up to a desired dimension. If in addition we
  add states in the form of letting the $\Omega_i$'s take vector
  spaces (Hilbert spaces, or some other algebraic structure) as
  values we enter the situation of topological quantum field theory
  which we will not pursue here, see Section \ref{sec:multi}.\\
\item[f)] \emph{Limits}.  In category theory we form limits and
  colimits of a collection of objects --- more precisely, given a
  functor
  \begin{equation*}
    F \colon I \to \catC
  \end{equation*}
  we form the colimit:
  \begin{equation*}
    \colim F = \colim c_i, \quad F(i) = c_i.
  \end{equation*}

  The colimit binds the collection or pattern of objects $c_i$ into
  one simple object reflecting the complexity of the pattern.  In this
  sense it is a bond in the hyperstructure framework if we drop the
  condition giving rise to the ``boundary'' maps $\partial_i$.

  If we require the $\partial_i$'s to exist, then the bond knows which
  objects it binds.  In the colimit this is not the case.  Hence we
  consider hyperstructures with and without $\partial_i$'s.

  Colimits may also be iterated.  For example we may consider
  situations where each $c_i$ already is a colimit of other colimits
  etc.  Expressed in a different way we consider a multivariable
  functor
  \begin{equation*}
    F \colon I_1 \times I_2 \times \cdots \times I_k \to \catC
  \end{equation*}
  and form the iterated colimit
  \begin{equation*}
    \calH \colon \quad \colim_{I_k}\ \cdots\ \colim_{I_2} \colim_{I_1}
    F.
  \end{equation*}
  This is clearly an iterated binding structure in the hyperstructure
  framework which we will discuss in the next section.  The colimits
  over the various indexing categories \ represent \ the \ bonds \ of
  \ the various \ levels, see \cite{EV,BEV} for a general context.
  For a categorical discussion of hyperstructure, see Appendix B in
  \citeasnoun{HAM}.
\end{itemize}

\section{A metaphor}
\label{sec:meta}
Let us illustrate using a metaphor what we mean by putting a
hyperstructure on an already existing structure, system or situation.

Suppose we are given a society or organization of agents, and we want
to act upon it in the manner of wielding political power, governing a
society or nation.  A possible procedure is: create a kind of
``political party'' organization.  A structural design of the
organization is needed, rules of action (``ideology'') and incentives
(``goals'').  The fundamental task is to create an organization --- a
``party'' --- starting with ``convinced'' individuals, then suitable
groups of individuals, groups of groups,$\ldots$

Basically this is putting a hyperstructure on the society of agents
which may act as an ideological amplifier from individuals to the
society.  This can be done independently of an existing societal
organization that one wants to act upon.  In such a hyperstructure the
bonds may depend on a goal (ideology) and incentives like solving
common problems, infrastructure, healthcare, poverty, etc.  In physics
it could be minimizing or releasing energy to obtain stability.  All
such factors will play a role in the build up of the hyperstructure in
the form of choosing bonds, states, etc.\ such that they support the
goals or ``ideology''.

Having established a hyperstructure, then let it act by the
``ideology'' in the sense that it should be instructible --- like a
superdemocracy.  Hence it may be instructed to maintain, replace or
improve the existing structure of society or change it to achieve
certain goals.  This is what political parties and other organizations
do.

A hyperstructure on a society (or space) will facilitate the
achievement of desired dynamics for agents or other objects through
the bonds which may act dynamically --- like fusion of old bonds to
new bonds.

If the hyperstructure is given as
\begin{equation*}
  \calH = \{B_0,B_1,\ldots,B_n\},
\end{equation*}
then one may initiate a dynamics at the lowest level
\begin{equation*}
  B_0 \to B_0'
\end{equation*}
which may require relatively few resources or little energy.  Then
these changes of actions and states will propagate through the higher
bonds
\begin{equation*}
  \begin{array}{c@{\ }c@{\ }c}
    B_1 & \to & B_1'\\
    & \vdots &\\
    B_n & \to & B_n'
  \end{array}
\end{equation*}
leading to a major desired change at the top level, depending on the
nature of $\calH$.  This is how a political organization works.  These
aspects are elaborated in Section \ref{sec:multi}.

The degree of detail of the hyperstructure will depend on the
situation and information available --- like how rich mathematically
the hyperstructure can be.  One may also ask:

What is the sociology of a space?  $\calH(X)$ represents the
``organization of $X$'' or ``the society of X''.  What can be obtained
by a political structure on $X$ depends on $X$ and
$\calH_{\text{political}}(X)$.

This shows through the metaphor how general the idea is and how it may
usefully apply in many situations.  Another interesting idea is how to
use this metaphor in the study of the genome as a society of genes.
One would like to put on a hyperstructure whose ideology should be to
maintain the structure (homeostasis), avoiding and discarding unwanted
growth like cancer.

How could one possibly create and represent such a ``genomic political
party''?  Possibly by an organized collection (hyperstructure) of
drugs (or external fields) that acts on bonds of genes.  The protein
P53 may already have such a role at a high level in an existing
hyperstructure.

\section{Binding structures}
\label{sec:binding}
We will now discuss the main issue in this paper, namely how to
organize a collection of objects (possibly with an existing structure)
into a new structure.

Let us assume that we are given a basic collection of objects
\begin{equation*}
  Z = \{z_i\}_{i\in I}
\end{equation*}
$I$ finite, countable or uncountable.  $Z$ may be the elements of a
set or a space.  Let us also assume that we are given a hyperstructure
of order $n$ in the sense of the previous section.
\begin{align*}
  \calH \colon \quad & \{X_0,X_1,\ldots,X_n\}\\
  & \{B_0,B_1,\ldots,B_{n - 1}\}
\end{align*}
In order to simplify the notation we do not write the $\Omega$'s and
$\partial$'s.\\

How to put an $\calH$-structure on $Z$?\\

This means describing how to bind the objects in $Z$ together into new
higher order objects following the pattern given by the bond structure
in $\calH$.  The idea is as follows.\\

We represent the collection $Z$ as a collection of elements in $X_0$
--- the basic set on which $\calH$ is built.

\begin{figure}[H]
  \centering
  \begin{tikzpicture}
    \node (zi) at (0,0){$\bullet$};
    \node at (0,-0.5){$z_i$};
    \node at (0,-1.75){$Z$};
    \node (zihat) at (5,0){$\bullet$};
    \node at (5,-0.5){$\hat{z}_i$};
    \node at (5,-1.75){$X_0$};

    \begin{scope}[scale=0.2,xshift=-4cm]
      \draw (0.811657,6.5691524) .. controls (0.104829654,5.8617663) and
        (0.0,-4.1397247) .. (0.61165696,-4.9308476).. controls
        (1.2233139,-5.7219706) and (7.5121303,-7.286697)
        .. (8.391657,-6.8108478) .. controls (9.271184,-6.334998) and
        (10.63916,-2.4085267) .. (10.291657,-1.4708476) .. controls
        (9.944154,-0.5331687) and (6.52478,-0.1069224)
        .. (6.231657,0.8491523) .. controls (5.9385343,1.805227) and
        (8.929111,4.5499115) .. (8.491657,5.4491525) .. controls
        (8.054203,6.3483934) and (4.007128,6.5140862)
        .. (3.011657,6.6091523) .. controls (2.016186,6.7042184) and
        (1.5184844,7.2765384) .. (0.811657,6.5691524);
    \end{scope}

    \begin{scope}[scale=0.2,xshift=21cm]
      \draw (0.811657,6.5691524) .. controls (0.104829654,5.8617663) and
        (0.0,-4.1397247) .. (0.61165696,-4.9308476).. controls
        (1.2233139,-5.7219706) and (7.5121303,-7.286697)
        .. (8.391657,-6.8108478) .. controls (9.271184,-6.334998) and
        (10.63916,-2.4085267) .. (10.291657,-1.4708476) .. controls
        (9.944154,-0.5331687) and (6.52478,-0.1069224)
        .. (6.231657,0.8491523) .. controls (5.9385343,1.805227) and
        (8.929111,4.5499115) .. (8.491657,5.4491525) .. controls
        (8.054203,6.3483934) and (4.007128,6.5140862)
        .. (3.011657,6.6091523) .. controls (2.016186,6.7042184) and
        (1.5184844,7.2765384) .. (0.811657,6.5691524);
    \end{scope}

    \path[->] (zi.north east) edge[bend right=-30]
    node[above]{$\widehat{\phantom{-}}$} (zihat.north west);
  \end{tikzpicture}
  \caption{}
  \label{fig:zx0}
\end{figure}

Hence we get a new collection of objects (or elements) in $X_0$ ---
\begin{equation*}
  \hat{Z} = \{\hat{z}_i\}_{i\in I}\subset X_0.
\end{equation*}
This is similar to the correspondence or analogy in Examples
\ref{ex:links} and \ref{ex:many} where particles are represented as
rings.  More on this later.

On $X_0$ we have bonds $B_0$ and can apply these to the new collection
$\hat{Z}$ in $X_0$.  Therefore we put a bond structure on $Z$ as
follows.

\begin{definition}
  \label{def:b0}
  $B_0(Z) = B_0(\hat{Z})$
\end{definition}

This means that we pull back the bonds from the hyperstructure $\calH$
on $X_0$ to $Z$.  Since $\hat{Z}\subset X_0$ we get an induced
hyperstructure on $\hat{Z}$.  This means that
\begin{equation*}
  \begin{array}{c}
    B_1(\hat{Z})\\
    B_2(\hat{Z})\\
    \vdots\\
    B_{n - 1}(\hat{Z})
  \end{array}
\end{equation*}
with the $\Omega_i$, $\partial_i$, $\Gamma_i$ coming along.

\begin{remark}
  If we already have a good hyperstructure on $Z$, we just keep it via
  the identity representation $(Z = X_0)$ and use it in the binding
  process we will describe.
\end{remark}

With a hyperstructure on the given collection $Z$ we can introduce new
higher order clusters and patterns of interactions.

\begin{definition}
  \label{def:clustering}
  An $\calH$-binding (or clustering) structure on $Z$ --- denoted by
  $\calH(Z)$ --- is given as follows:

  Let $S_0\subset Z$ and $b_0\in B_0(S_0)$, $z\in Z$ is an element of
  a $b_0$-cluster ($\clu(b_0)$) if $z\in S_0$.

  Furthermore, let $S_1\subset B_0(Z)$ and $b_1\in B_1(S_1)$, then
  $b_0\in B_0(Z)$ is an element of a $b_1$-cluster ($\clu(b_1)$) if
  $b_0\in S_1$.

  If $z\in \clu(b_0)$, and $b_0\in \clu(b_1)$, then we have a second
  order clustering and write $z\in \clu(b_0,b_1)$.  In the same way we
  proceed to $n$-th order clustering by requiring
  \begin{equation*}
    z\in \clu(b_0,b_1,\ldots,b_{n - 1})
  \end{equation*}
  in an obvious extension of the notation.
\end{definition}

This describes the general $\calH$-binding (or clustering) principle.
The same principle applies to the extended representation picture:
\begin{center}
  \begin{tikzpicture}
    \node (Zup) at (0,1.5){$Z$};
    \node (X0) at (1.5,1.5){$X_0$};
    \node (Zupr) at (3,1.5){$Z$};

    \node (Zdown) at (0,0){$Z$};
    \node (Xn) at (1.5,0){$X_n$};
    \node (Zdownr) at (3,0){$Z$};

    \begin{scope}[font=\scriptsize,midway,above]
      \draw[->] (Zup) -- node{$R_0$} (X0);
      \draw[->] (X0) -- node{$I_0$} (Zupr);
      \draw[->] (Zdown) -- node{$R_n$} (Xn);
      \draw[->] (Xn) -- node{$I_n$} (Zdownr);
      \draw[dotted] (X0) -- (Xn);
    \end{scope}
  \end{tikzpicture}
\end{center}
interpreted in the natural way: $R$ giving a binding structure and $I$
inducing a ``parametrization'' or decomposition by taking inverse
images.  The figure indicates that we may represent or induce at any
level, but most of the time we use level $0$.  One may construct a
decomposition of $Z$ via a hyperstructure $\calH$, by starting with
the top bonds $B_n$ (reverse the direction of the binding process).
Furthermore, one may then bind the lowest level (smallest) elements of
the decomposition ($B_0$-bonds) to a new hyperstructure
$\hat{\calH}$.  The situation may also be extended to the $R$'s and
$I$'s being of relational character.

The $\calH$-binding principle that we have described is in a way also
a Transfer Principle of Organization --- showing how to transfer
structure and organization from one universe to another (this is more
general than functors between categories).  For example one may use it
to transfer deep geometrical bonds to other interacting systems, like
particle systems as described in \citeasnoun{NS} and
\citeasnoun{manybody}.

The idea may be easier to grasp in the case that the hyperstructure is
given by a composition:
\begin{equation*}
  \calH \colon \quad S_1 \xleftarrow{\varphi_1} S_2
  \xleftarrow{\varphi_2} \cdots \xleftarrow{\varphi_{n - 1}} S_n.
\end{equation*}
In this case a given collection $Z$ should be represented in $S_n$.
For simplicity let us consider the identity representation $Z = S_n$.
The elements of $S_{n - 1}$ represent $B_0$, hence $z\in \clu(s_{n -
  1})$ if $z\in \varphi_{n - 1}^{-1} (s_{n - 1})$.

Similarly
\begin{equation*}
  s_{n - 1}\in \clu(s_{n - 2}) \quad \text{if} \quad s_{n - 1}\in
  \varphi_{n - 2}^{-1} (s_{n - 2})
\end{equation*}
and $2$nd order clusters are formed $\clu(s_{n - 1},s_{n - 2})$.
Hence $z\in \clu(s_{n - 1},s_{n - 2},\ldots,s_1)$ gives an $n$-th
order clustering scheme of $Z$.

This discussion shows what we mean by putting a hyperstructure on a
collection of objects.  In a way the hyperstructure acts as a
parametrization of the new objects formed by higher bonds.  An
important point here is the choice of representations, and how to
choose them in an interesting and relevant way.  $Z$ may be any
collection of elements, subsets, subspaces, elements of a
decomposition, etc.  $\calH(Z)$ organizes $Z$ into a ``society'' of
objects.

Putting a hyperstructure on a ``situation'' is not meant to be
restricted to a set or a space, but could be a non-mathematical
``situation'', a category of some kind or in general an already
existing hyperstructure (or families of such).  The transfer of
binding structures may be considered as a kind of generalized
``functoriality''.

In fact one may think of a given situation $S$ in representing level
$0$, and then always look for a higher order associated situation in
the form of a hyperstructure --- $\calH(S)$, which may give a better
understanding of the given situation (or object).  For example:
\begin{align*}
  &S \colon \text{a set $X_0$} \rightsquigarrow \calH(X_0)\\
  &S \colon \text{a category $\catC_0$} \rightsquigarrow
  \calH(\catC_0)
\end{align*}
This idea is a bit similar to the idea of associating complexes and
resolutions to groups, modules, algebras, categories, etc.\ ---
derived objects (or situations).

Sometimes, a representation is given in a natural way by the nature of
the object.  For example molecules being represented by their own
geometric form.  On the other hand one may choose more abstract
representations embedding the collection in a universe rich with
structures and possibilities for interesting bindings and
interactions.

In this way we put a hyperstructure on a situation by somehow inducing
one from a known one.  If a good model does not exist, one may have to
create a new suitable hyperstructure for further use.  In either case
the hyperstructure enables us to form organized Abstract Matter in the
sense of \citeasnoun{HAM} in order to handle a collection of objects or a
situation and achieve certain goals.  Binding structures represent the
\emph{Principle of how we make Things}.

As pointed out in \citeasnoun{NS} one may study the collection in a
selected universe and then ask whether the ``abstract'' binding
structure may be realized in the original universe.  Hence we may use
an $\calH$-structure to synthesize new bond states.

This is exactly the situation we study in \citeasnoun{NS} where $Z$ is
a collection of particles in a cold gas.  We represent the particles
by rings as in Examples \ref{ex:links} and \ref{ex:many}.  The
hyperstructure of higher order links --- $\calH(\text{Rings})$ ---
pulls back to a hyperstructure of the particle collection ---
$\calH(\text{Particles})$.  It is a verified fact, see \citeasnoun{NS}
and references therein, that to Brunnian (or Borromean) rings there
corresponds quantum mechanical states --- Efimov states with the same
binding patterns.

On this basis it is tempting to suggest that there is a similar
correspondence between higher order links and states corresponding to
the higher order clusters of bound particles --- higher order Brunnian
states.  This is the main suggestion of \citeasnoun{NS}.

For example the binding pattern of a second order Brunnian ring
$2B(3,3)$ suggests a bond or particle state of $3$ particles bound to
trimers and $3$ trimers bound to a single state.  This is the key idea
in our guiding examples of the general setting of binding structures
we have introduced.

We may summarize our discussion in the following figure:

\begin{figure}[H]
  \centering
  \begin{tikzpicture}
    \begin{scope}
      \draw (0,0) rectangle (2.5,5);
      \foreach \y in {1,2,4}{
        \draw (0,\y) -- (2.5,\y);
      }

      \node at (1.25,-0.6){$\calH_{R_0}(Z)$};
      \node at (1.25,0.4){$\calH_{R_0}^0(Z) = Z$};
      \node at (1.25,1.4){$\calH_{R_0}^1(Z)$};
      \node at (1.25,3){$\vdots$};
      \node at (1.25,4.4){$\calH_{R_0}^n(Z)$};

      \draw[->] (3,0.4) -- node[midway,above]{$R_0$} (4.5,0.4);
      
      \foreach \y in {1.5,4.5}{
        \draw[dashed,->] (3,\y) -- (4.5,\y);
      }
    \end{scope}

    \begin{scope}[xshift=5cm]
      \draw (0,0) rectangle (2.5,5);
      \foreach \y in {1,2,4}{
        \draw (0,\y) -- (2.5,\y);
      }

      \node at (1.25,-0.6){$\calH(X_0)$};
      \node at (1.25,0.4){$\calH^0(X_0) = X_0$};
      \node at (1.25,1.4){$\calH^1(X_0) = X_1$};
      \node at (1.25,3){$\vdots$};
      \node at (1.25,4.4){$\calH^n(X_0) = X_n$};

      \draw[->] (3,0.4) -- node[midway,above]{$I_0$} (4.5,0.4);
      
      \foreach \y in {1.5,4.5}{
        \draw[dashed,->] (3,\y) -- (4.5,\y);
      }
    \end{scope}

    \begin{scope}[xshift=10cm]
      \draw (0,0) rectangle (2.5,5);
      \foreach \y in {1,2,4}{
        \draw (0,\y) -- (2.5,\y);
      }

      \node at (1.25,-0.6){$\calH_{I_0}(Z)$};
      \node at (1.25,0.4){$\calH_{I_0}^0(Z) = Z$};
      \node at (1.25,1.4){$\calH_{I_0}^1(Z)$};
      \node at (1.25,3){$\vdots$};
      \node at (1.25,4.4){$\calH_{I_0}^n(Z)$};
    \end{scope}
  \end{tikzpicture}
  \caption{}
  \label{fig:basic_prin}
\end{figure}

\noindent illustrating the two basic principles:\\

\begin{itemize}
\item[(I)] The Hyperstructure Principle --- which is an organizing
  principle and a guiding principle for structural architecture and
  engineering.\\
\item[(II)] The Transfer Principle --- which is a way to transfer a
  hyperstructure.\\
\end{itemize}

For the situation in Example \ref{ex:many} the figure looks like:

\begin{figure}[H]
  \centering
  \begin{tikzpicture}[scale=0.7]
    \begin{scope}
      \draw (0,0) rectangle (5,9);
      \foreach \y in {3,6}{
        \draw (0,\y) -- (5,\y);
      }

      \begin{scope}[xshift=2.5cm,yshift=7.5cm,scale=0.55]
        \draw[very thick] (0,0) circle(2.5cm);
        
        \begin{scope}[scale=0.3,xshift=-4.25cm,yshift=2.5cm]
          \draw[very thick] (0,0) circle(2.5cm);
          \draw[very thick, blue] (-0.25,-1.6) circle(0.5cm);
          \draw[very thick, red] (-1.1,1.2) circle(0.5cm);
          \draw[very thick, green] (1.4,0.8) circle(0.5cm);
        \end{scope}
        
        \begin{scope}[scale=0.3,xshift=4.5cm,yshift=1.2cm]
          \draw[very thick] (0,0) circle(2.5cm);
          \draw[very thick, blue] (-0.25,-1.6) circle(0.5cm);
          \draw[very thick, red] (-1.1,1.2) circle(0.5cm);
          \draw[very thick, green] (1.4,0.8) circle(0.5cm);
        \end{scope}
        
        \begin{scope}[scale=0.3,xshift=-0.5cm,yshift=-4cm]
          \draw[very thick] (0,0) circle(2.5cm);
          \draw[very thick, blue] (-0.25,-1.6) circle(0.5cm);
          \draw[very thick, red] (-1.1,1.2) circle(0.5cm);
          \draw[very thick, green] (1.4,0.8) circle(0.5cm);
        \end{scope}
      \end{scope}

      \begin{scope}[xshift=2.5cm,yshift=4.5cm,scale=0.65]
        \begin{scope}[scale=0.3,xshift=-4.25cm,yshift=2.5cm]
          \draw[very thick] (0,0) circle(2.5cm);
          \draw[very thick, blue] (-0.25,-1.6) circle(0.5cm);
          \draw[very thick, red] (-1.1,1.2) circle(0.5cm);
          \draw[very thick, green] (1.4,0.8) circle(0.5cm);
        \end{scope}
        
        \begin{scope}[scale=0.3,xshift=4.5cm,yshift=1.2cm]
          \draw[very thick] (0,0) circle(2.5cm);
          \draw[very thick, blue] (-0.25,-1.6) circle(0.5cm);
          \draw[very thick, red] (-1.1,1.2) circle(0.5cm);
          \draw[very thick, green] (1.4,0.8) circle(0.5cm);
        \end{scope}
        
        \begin{scope}[scale=0.3,xshift=-0.5cm,yshift=-4cm]
          \draw[very thick] (0,0) circle(2.5cm);
          \draw[very thick, blue] (-0.25,-1.6) circle(0.5cm);
          \draw[very thick, red] (-1.1,1.2) circle(0.5cm);
          \draw[very thick, green] (1.4,0.8) circle(0.5cm);
        \end{scope}
      \end{scope}

      \foreach \x in {1.25,2.5,3.75}{
        \draw[thick,red] (\x,2.5) circle(0.3cm);
        \draw[thick,green] (\x,1.5) circle(0.3cm);
        \draw[thick,blue] (\x,0.5) circle(0.3cm);
      }

      \node at (2.5,9.8){State hyperstructure:};
      \node at (2.5,-0.6){$Z = \text{$9$ particle states}$};

      \foreach \y in {4.5,7.5}{
        \draw[line width=1.5,->] (6,\y) -- (9,\y);
      }

      \draw[line width=1.5,->] (6,1.5) -- node[midway,below]{$R_0$}
        (9,1.5);
    \end{scope}

    \begin{scope}[xshift=10cm]
      \draw (0,0) rectangle (5,9);
      \foreach \y in {3,6}{
        \draw (0,\y) -- (5,\y);
      }

      \begin{scope}[knot/.style={thin knot},scale=0.25]
        \begin{scope}[xshift=5cm,yshift=20cm]
          \setbrstep{.1}
          \brunnian{2.5}{3}
        \end{scope}

        \begin{scope}[xshift=15cm,yshift=20cm]
          \setbrstep{.1}
          \brunnian{2.5}{3}
        \end{scope}

        \begin{scope}[xshift=10cm,yshift=16cm]
          \setbrstep{.1}
          \brunnian{2.5}{3}
        \end{scope}
      \end{scope}

      \foreach \x in {1.25,2.5,3.75}{
        \draw[thick,red] (\x,2.5) circle(0.3cm);
        \draw[thick,green] (\x,1.5) circle(0.3cm);
        \draw[thick,blue] (\x,0.5) circle(0.3cm);
      }

      \node at (2.5,9.8){Link hyperstructure:};
      \node at (2.5,-0.6){$X_0 = \text{$9$ rings}$};
    \end{scope}
  \end{tikzpicture}
  \raisebox{5.25cm}[0pt][0pt]{
    \hspace*{7cm}
    \includegraphics[scale=0.09]{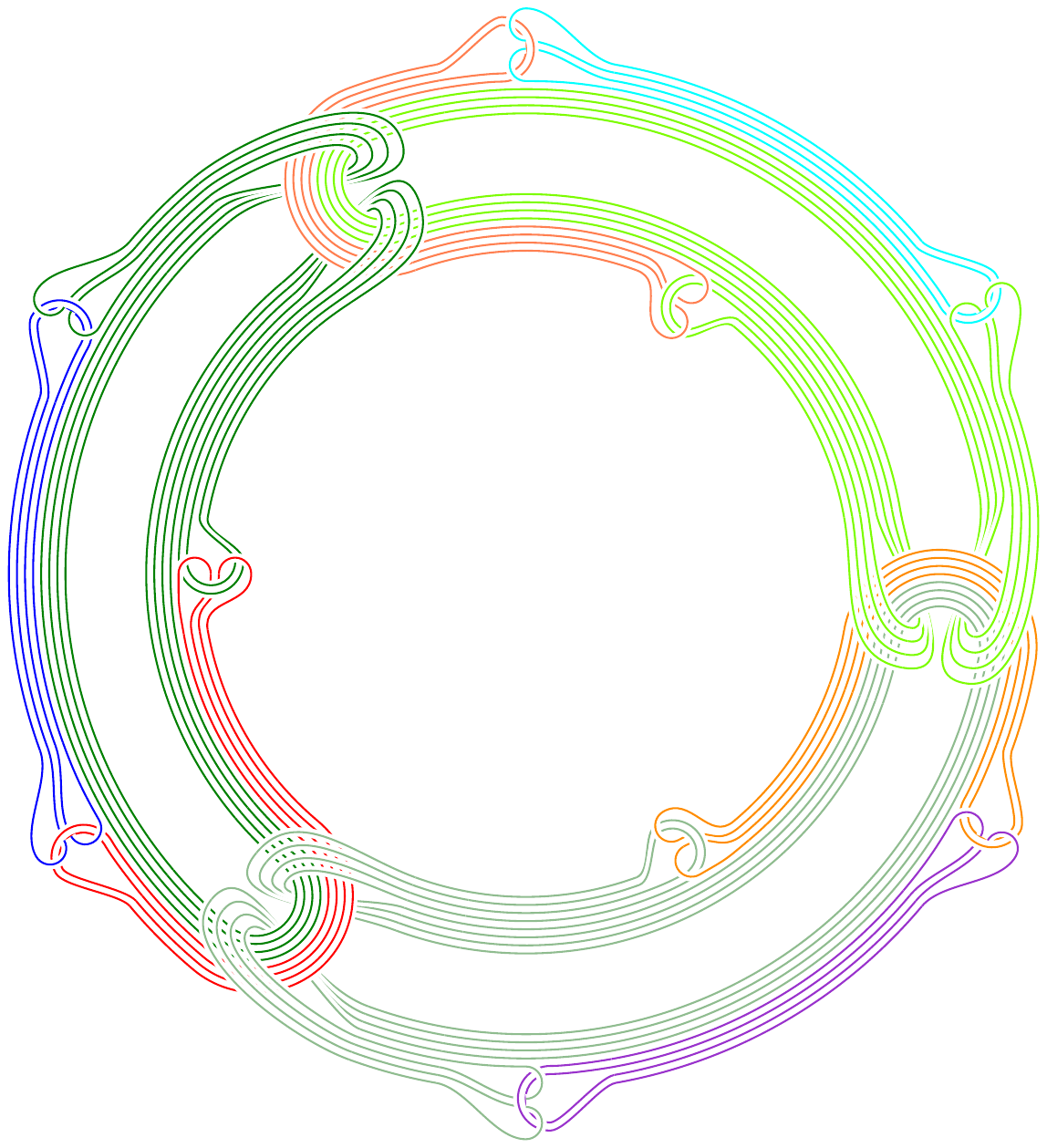}
  }
  \caption{}
  \label{fig:example2}
\end{figure}

The binding structures and diagrams in Figure \ref{fig:basic_prin}
described here may be considered as extensions of pasting diagrams in
higher categories.  A sheaf type formulation has been given in
Appendix B in \citeasnoun{HAM}.  Using the Principles described here
one may induce an action on a totality $Z$, by acting on individual
elements and letting the action propagate through $\calH(Z)$ to the
top level as in political processes and in social and business
organizations.  In this way many small actions may lead to major
global actions and change of state.  Hence $\calH$ acts as an
amplifier.

\section{Analysis and synthesis}
\label{sec:analsyn}
How do we synthesize new objects or structures from old ones?  A
common procedure both in nature and science is to bind objects
together to form new objects with new properties, then use these
properties in forming new bonds and new higher order objects.  This is
precisely what a hyperstructure does!

Let $Z$ be a collection of objects.  By putting a hyperstructure on
$Z$ --- $\calH(Z)$ we have a binding scheme of new higher order
synthetic objects.  If $\calH(Z)$ is pulled back from a given
structure, the problem may be to tune the environment of $Z$ in such a
way that the binding pattern of $\calH(Z)$ may be realized.

On the other hand $Z$ may be considered as a global object that we
want to analyze by decomposing it into smaller and smaller pieces.  By
putting a hyperstructure on $Z$ --- $\calH(Z)$ --- we have seen in the
previous section how it gives rise to a higher order clustering
decomposition of $Z$.

It is interesting to notice that if we take the smallest pieces in the
decomposition (lowest level elements) as our basic set $Z(\dec)$, we
may put a new hyperstructure $\widehat{\calH}$ on it and recapture $Z$
as the top level of $\widehat{\calH}\left(Z(\dec)\right)$.

This shows that hyperstructures are useful in the \emph{synthesis} of
new collections of objects from given ones and in the \emph{analysis}
of them as well.  It is very useful to put a hyperstructure on a
collection of objects in order to manipulate the collection towards
certain goals.  We will discuss various applications in the following.

\section{Applications of the binding structure}
\label{sec:app}
Many interactions in science and nature may be described and handled
as organized and structured collections of objects --- in certain
contexts called many-body systems.  In the following we would like to
point out that hyperstructures of bindings may be interesting and
useful in many areas.

\begin{itemize}
\item[a)] \emph{Physics}.  $\calH$-structures of many body systems
  (particles) may --- as we have already discussed --- give rise to
  new and exotic states of matter ($\calH$-states), see Examples
  \ref{ex:links} and \ref{ex:many}.\\
\item[b)] \emph{Chemistry}.  Hyperstructures like higher order links
  are interesting models for the synthesis of new molecules and
  materials --- like higher order Brunnian rings, see
  \citeasnoun{NSCS} and \citeasnoun{BS}.

  More generally we may consider a hyperstructure $\calH$ where the
  bonds are spaces like manifolds or CW-complexes built up of cells.

  Let us think of a collection of molecules each represented by a
  point in space.  Then they may form

  \begin{figure}[H]
    \centering
    \begin{tikzpicture}[scale=0.8]
      \node at (0,0){Chains:};
      \node at (1.75,-0.75){[or a ring $\bigcirc$ if we prefer]};
      \node at (8,0){represented by};
      \node at (12.4,0){(an interval)};
      \node at (13.8,-0.75){a $1$ dim geometric object};

      \foreach \a in {1,2}{
        \draw (\a,0) -- (\a+1,0);
        \draw[fill=black] (\a,0) circle(0.075cm);
      }

      \draw[fill=black] (3,0) circle(0.075cm);
      \draw (3,0) -- (3.625,0);
      \node at (4,-0.05){$\cdots$};
      \draw (4.375,0) -- (6,0);
      \draw[fill=black] (5,0) circle(0.075cm);
      \draw[fill=black] (6,0) circle(0.075cm);

      \draw (9.9,0) -- (10.9,0);
    \end{tikzpicture}
    \caption{}
    \label{fig:chains}
  \end{figure}

  This representation increases the dimension by one.  Similarly:

  \begin{figure}[H]
    \begin{tikzpicture}[scale=1]
      \foreach \l in {0,0.5,1,1.5}{
        \draw (0,\l) -- (2,\l);
      }

      \node at (3.75,0.75){(bound together)};
      \node at (7.5,0.75){may be represented by};

      \draw (10,0) rectangle(12,1.5);

      \node at (13.5,0.75){(a rectangle)};
    \end{tikzpicture}
    \caption{}
    \label{fig:rectangle}
  \end{figure}

  Collection of rectangles forming new chains:

  \begin{figure}[H]
    \begin{tikzpicture}
      \draw (0,0) rectangle (2,2);

      \draw (0.5,2.5) -- (2.5,2.5) -- (2.5,0.5);
      \draw (0.5,2.5) -- (0.5,2.0675);
      \draw[dashed] (0.5,1.875) -- (0.5,0.5) -- (1.875,0.5);
      \draw (2.0675,0.5) -- (2.5,0.5);

      \draw (1,3) -- (3,3) -- (3,1);
      \draw (1,3) -- (1,2.5675);
      \draw[dashed] (1,2.375) -- (1,1) -- (2.375,1);
      \draw (2.5675,1) -- (3,1);
      \draw[dotted] (1,3) -- (0.5,2.5);
      \draw[dotted] (3,3) -- (2.5,2.5);
      \draw[dotted] (3,1) -- (2.5,0.5);

      \path (5,1.25) edge[->,bend right=-30]
        node[midway,above,yshift=0.25cm]{may be represented by $3$
          dim box} (9,1.25); 

      \draw (11,0) rectangle (13,2);
      \draw (11,2) -- (12,3) -- (14,3) -- (14,1) -- (13,0);
      \draw (13,2) -- (14,3);
      \draw[dashed] (11,0) -- (12,1) -- (14,1);
      \draw[dashed] (12,3) -- (12,1);
    \end{tikzpicture}
    \caption{}
    \label{fig:box}
  \end{figure}

  Then we form chains of $3$ dimensional boxes again to be represented
  by a $4$ dimensional box, etc.  We may also introduce holes and we
  may continue up to a desired dimension.

  Then we may glue the molecular cells together following the
  topological patterns (for example homotopy type) of the bond spaces
  $B_n$ in each dimension.

  In this way we get organized molecules in three dimensions with the
  structure induced from a higher dimensional binding structure in
  $\calH$.  Clearly many other similar representations are possible.
  This is a very useful and important principle.

  We could also have molecules representing the bonds as in the
  following figures:

  \begin{figure}[H]
    \centering
    \begin{tikzpicture}
      \begin{scope}
        \foreach \y in {0,0.75,1,,1.25,1.5}{
          \draw (0,\y) node[left,yshift=-0.05cm]{$\cdots$} -- (1,\y);
        }

        \node at (-0.375,0.45){$\vdots$};
        \node at (0.5,0.45){$\vdots$};
        \node[below,yshift=-0.15cm] at (0.5,0){$b_0$};
      \end{scope}

      \begin{scope}[xshift=4cm]
        \foreach \y in {0,0.75,1,,1.25,1.5}{
          \draw (0,\y) -- (1,\y);
        }
        
        \draw (0,0) -- (0,1.5);
        \node at (0.5,0.45){$\vdots$};
        \node[below,yshift=-0.15cm] at (0.5,0){$b_1$};
      \end{scope}

      \begin{scope}[xshift=8cm]
        \foreach \y in {0,0.75,1.125,1.5}{
          \draw (0,\y) -- (1,\y);
          \foreach \x in {0.1,0.35,0.6,0.85}{
            \draw (\x,\y) -- (\x+0.15,\y+0.15);
          }
        }

        \draw (0,0) -- (0,1.5);
        \node at (0.5,0.5){$\vdots$};
        \node[below,yshift=-0.15cm] at (0.5,0){$b_2$};
      \end{scope}

      \begin{scope}[xshift=12cm]
        \foreach \y in {0,0.75,1,1.25,1.5}{
          \draw[decorate,decoration=zigzag] (0,\y) -- (1,\y);
        }

        \draw (0,0) -- (0,1.5);

        \node at (0.5,0.5){$\vdots$};
        \node[below,yshift=-0.15cm] at (0.5,0){$b_n$};
      \end{scope}

      \foreach \x in {2,6,10}{
        \path (\x,1) edge[->,bend right=-30] (\x+1,1);
      }
    \end{tikzpicture}
    \caption{}
    \label{fig:molecules}
  \end{figure}

  \noindent describing the process:
  
  \begin{figure}[H]
    \centering
    \begin{tikzpicture}[text centered]
      \node at (0,1.25){$\calM$};
      \node[text width=2cm] at (0,0){collection of molecules};
      \node at (5,1.25){$\calH(\calM)$};
      \node[text width=4cm] at (5,0){an $\calH$-structured bond
        collection};

      \path (0.5,1.3) edge[->,bend right=-15] (4,1.3);
    \end{tikzpicture}
    \caption{}
    \label{fig:process}
  \end{figure}

\item[c)] \emph{Social and economic systems}.  We may here consider
  populations of individuals, social or economic units.  Then it may
  be useful to consider them as many body systems in the physical
  examples and introduce higher order binding structures, see examples
  in \citeasnoun{NS} and Section \ref{sec:meta}.

  One may for example discuss Brunnian investments of $n$ agents and
  continuation to higher order which may be interesting in certain
  contexts.\\
\item[d)] \emph{Biology}.  Here we may consider collections of genes,
  cells, pathways, neurons, etc.\ as many body systems and bind them
  together in new ways according to a given $\calH$-structure.  For
  example wihtin tissue engineering one may make $\calH$-type tissues
  for various purposes.\\
\item[e)] \emph{Logic}.  We may introduce $\calH$-type bindings of
  logical types and data-structures.  New ``laws of thought'' are
  possible based on a logic of $\calH$-type bindings as ``deductions''
  and states/observations in $\calH$ representing the semantics.\\
\item[f)] \emph{Networks}.  In \citeasnoun{NSCS} we argued that in many
  situations networks are inadequate and should be replaced by
  hyperstructures.  Pairwise binding or interactions would then be
  replaced by $\calH$-bindings.  Look at the Brunnian hyperstructure
  of links as in the introduction.\\
\item[g)] \emph{Brain systems}.  Extend natural and artificial neural
  networks to $\calH$-structures of neurons as follows.

  Let $Z$ be a collection of real or abstract neurons.  Then the
  $\calH(Z)$ binding structure represent new interaction patterns
  ``parametrized'' by $\calH$, possibly representing new types of
  higher order cognitive functions and properties.  See also
  \citeasnoun{Cog}.\\
\item[h)] \emph{Correlations}.  One may think of correlations as
  relations and bindings of variables.  An interesting possibility
  would be to extend pair correlations to $\calH$-type correlations of
  higher clusters.  They could possibly have Brunnian properties as
  follows:
  \begin{equation*}
    \corr(X,Y) = \corr(Y,Z) = \corr(X,Z) = 0, \quad \text{but} \quad
    \widehat{\corr}(X,Y,Z) \neq 0.
  \end{equation*}
  This is in analogy with cup products and Massey products in the
  study of Brunnian links.  To detect higher order Brunnian linking
  one introduces higher order Massey products.  In the correlation
  language this would mean putting
  \begin{equation*}
    \hat{A} = \widehat{\corr}(X_1,Y_1,Z_1) = 0, \quad \hat{B} =
    \widehat{\corr}(X_2,Y_2,Z_2) = 0, \quad \hat{C} =
    \widehat{\corr}(X_3,Y_3,Z_3) = 0
  \end{equation*}
  and finding a $\hat{\hat{C}}$ such that
  $\hat{\hat{C}}(\hat{A},\hat{B},\hat{C}) \neq 0$ represents a second
  order correlation.\\
\item[i)] \emph{Mathematics}.  As already indicated in some of the
  examples collections of mathematical objects may also be bound
  together in new ways modelled or para-metrized by a hyperstructure,
  for example collections of spaces like manifolds and cell-complexes
  along with gluing bonds.  

  Assume that we have a collection of spaces $\calL = \{L_i\}$
  organized in a well-defined hyperstructure $\calH(\calL)$.  If we
  have another collection of spaces or objects $\calX = \{X_i\}$ we
  may then induce an $\calH$-structure on $\calX$ --- $\calH(\calX)$
  and use it to study the collection $\calX$.

  \begin{figure}[H]
    \centering
    \begin{tikzpicture}
      \begin{scope}[scale=0.5,rotate=135,yshift=1.75cm]
        \draw (1.4652985,5.512821) .. controls (2.1881015,6.2038755)
          and (6.217555,5.638671) .. (6.9452987,4.9528213) .. controls
          (7.673042,4.2669716) and (7.8358297,2.6816413)
          .. (7.4652987,1.7528212) .. controls (7.0947676,0.8240011)
          and (5.517234,0.9829246) .. (5.2052984,0.032821205)
          .. controls (4.8933635,-0.91728216) and
          (6.6252465,-3.2035027) .. (6.4452987,-4.1871786) .. controls
          (6.2653503,-5.170855) and (5.6661596,-6.210148)
          .. (4.7052984,-6.487179) .. controls (3.7444375,-6.7642097)
          and (3.4615378,-5.8638744) .. (2.4852986,-5.6471786)
          .. controls (1.5090593,-5.4304833) and
          (1.2496852,-6.2882943) .. (0.62529856,-5.507179) .. controls
          (0.0,-4.7260633) and (1.6402918,0.6538345)
          .. (1.6852986,1.6528212) .. controls (1.7303053,2.6518078)
          and (0.7424957,4.821767) .. (1.4652985,5.512821);
      \end{scope}

      \begin{scope}[xshift=-3.25cm,yshift=-1.25cm,scale=0.65]
        \draw (0,0) .. controls (1,0.5) and (-0.5,1.25) .. (0,2);
      \end{scope}

      \begin{scope}[xshift=-2.25cm,yshift=-0.5cm,scale=0.65]
        \draw (0,0) .. controls (1,0.5) and (-0.5,1.25) .. (0,2);
      \end{scope}

      \node at (-1.5,0.5){$X_i$};

      \begin{scope}[xshift=-0.25cm,yshift=1cm,scale=0.65]
        \draw (0,0) .. controls (1,0.5) and (-0.5,1.25) .. (0,2);
      \end{scope}

      \begin{scope}[scale=0.5,xshift=-11cm]
        \draw (27.885298,5.972821) .. controls (27.2742,6.764375) and
          (20.968887,6.534102) .. (20.185299,5.9128213) .. controls
          (19.40171,5.29154) and (19.228222,4.2212214)
          .. (19.545298,3.2728212) .. controls (19.862375,2.324421)
          and (21.11391,2.886598) .. (21.225298,1.8928212) .. controls
          (21.336687,0.8990443) and (16.50174,-3.090676)
          .. (16.585299,-4.0871787) .. controls (16.668858,-5.0836816)
          and (18.332634,-6.486274) .. (19.325298,-6.6071787)
          .. controls (20.317963,-6.7280836) and
          (24.766945,-4.0730786) .. (25.385298,-3.2871788) .. controls
          (26.003653,-2.501279) and (28.496397,5.1812673)
          .. (27.885298,5.972821);
      \end{scope}

      \begin{scope}[xshift=4cm,yshift=-2.75cm,scale=0.65,rotate=10]
        \draw (0,0) .. controls (1,0.5) and (-0.5,1.25) .. (0,2);
      \end{scope}

      \begin{scope}[xshift=5cm,yshift=-1.75cm,scale=0.65,rotate=10]
        \draw (0,0) .. controls (1,0.5) and (-0.5,1.25) .. (0,2);
      \end{scope}

      \node at (5.5,-0.75){$L_i$};

      \begin{scope}[xshift=6.5cm,yshift=0.25cm,scale=0.65,rotate=10]
        \draw (0,0) .. controls (1,0.5) and (-0.5,1.25) .. (0,2);
      \end{scope}

      \node at (0,4){$\bond_\calX$};
      \node at (8,4){$\bond_\calL$};
    \end{tikzpicture}
    \caption{}
    \label{fig:calL}
  \end{figure}
  Interesting properties for $\calH(\calL)$ may be asked for
  $\calH(\calX)$ like Brunnian properties in a categorical setting,
  see \citeasnoun{NS}.

  We may take a family of spaces, for example simplicial complexes:
  \begin{equation*}
    \calV = \{V_i\}
  \end{equation*}
  and represent them in a family of manifolds
  \begin{equation*}
    \calM = \{M_i\}
  \end{equation*}
  on which there is a hyperstructure $\calH(\calM)$ with given bonds
  such that
  \begin{equation*}
    V_i \to M_i
  \end{equation*}
  and
  \begin{figure}[H]
    \centering
    \begin{tikzpicture}
      \begin{scope}[scale=0.8]
        \draw (0.15224439,-1.3088658) .. controls
        (-0.001173694,-2.297027) and (0.094758034,-2.0608437)
        .. (1.0322444,-2.4088657) .. controls (1.9697307,-2.7568877)
        and (11.833596,-1.1438532) .. (12.352244,-0.28886575)
        .. controls (12.8708935,0.56612176) and (12.43143,1.6672896)
        .. (11.592244,2.2111342) .. controls (10.753058,2.754979) and
        (4.9292545,1.6004366) .. (3.9922445,1.2511343) .. controls
        (3.0552344,0.90183187) and (0.30566248,-0.32070437)
        .. (0.15224439,-1.3088658);
      \end{scope}

      \begin{scope}[xshift=1.25cm,yshift=-1.75cm,scale=0.65,rotate=5]
        \draw (0,0) .. controls (1,0.5) and (-0.5,1.25) .. (0,2);
      \end{scope}

      \begin{scope}[xshift=2.75cm,yshift=-1.375cm,scale=0.65,rotate=5]
        \draw (0,0) .. controls (1,0.5) and (-0.5,1.25) .. (0,2);
      \end{scope}      

      \begin{scope}[xshift=4.25cm,yshift=-1cm,scale=0.65,rotate=5]
        \draw (0,0) .. controls (1,0.5) and (-0.5,1.25) .. (0,2);
      \end{scope}

      \begin{scope}[xshift=7.75cm,yshift=-0.3cm,scale=0.65,rotate=5]
        \draw (0,0) .. controls (1,0.5) and (-0.5,1.25) .. (0,2);
      \end{scope}

      \node at (5,-0.25){$M_i$};
      \node at (10,2){$B_0$};
    \end{tikzpicture}
    \caption{}
    \label{fig:manifold}
  \end{figure}
  
  \noindent $B_0$ being a bond for $\{M_i\}$, for example by $M_i
  \subset B_0$.  Then we may say that the $\{V_i\}$ binds by a
  pullback bond $\hat{B}_0$.  Similar for higher bonds.  In this way
  one may introduce geometric structures for each level which
  otherwise may have been difficult.  In the sense that
  hyperstructures extend the various notions of higher categories, one
  may also introduce $\calH$-type bundles and stacks with transition
  and gluing morphisms replaced by appropriate bonds.  Then one may
  hope to extend bundle notions like connections, curvature and
  holonomy in suitable contexts.
\end{itemize}

Let us follow up the examples and discussion of hyperstructures in
Section \ref{sec:hyper}.\\

How to produce hyperstructures?\\

We have seen that compositions of maps
\begin{equation*}
  S_1 \leftarrow S_2 \leftarrow \cdots \leftarrow S_n
\end{equation*}
naturally lead to hyperstructures on $S_n$.  This may also be extended
to a situation of compositions of functors and sub-categories.

Geometric bonds are basically constructed by binding families of
spaces or sub-spaces of an ambient space:
\begin{equation*}
  \left\{A\phantom{^\ast}, \{A(\omega)\}\right\}, \quad \omega =
  (i_1,i_2,\ldots,i_n),
  \quad i_j\in I_j
\end{equation*}
where
\begin{equation*}
  A \supset A(i_1) \supset A(i_1,i_2) \supset \cdots \supset
  A(\omega).
\end{equation*}
are successive bonds.  Here the $A(\omega)$'s may be general spaces,
manifolds, etc.

Generalized link and knot theory may be viewed as the study of
embeddings of topological spaces in other topological spaces, but
hyperstructures encompass much more.  Still for geometric
hyperstructures one may consider using and extending the mathematical
theory of links and knots (quantum versions, quandles, etc., see
\citeasnoun{N}) to the study of geometric hyperstructures.

Manifolds with singularities as introduced in \citeasnoun{Sing} are
also represented by such bond systems $A = \{A(\omega)\}$.  So is also
the Brunnian link hyperstructure
\begin{equation*}
  B = \{B(\omega)\} = \{B(n_1,n_2,\ldots,n_k)\}.
\end{equation*}
This applies to structures in general --- for example algebraic,
topological, geometric, logical and categorical --- presented as
follows
\begin{equation*}
  \catC = \{\catC(\omega)\},
\end{equation*}
where structure $\catC(\omega)$ binds the structures $\catC(\omega,i)$
for example as substructures --- like in higher order links and
many-body systems.

This may be viewed as a kind of many-body problem of general
structures, and represent a simple organizing principle for them.

All this shows that there is a plethora of possible applications of
hyperstructured binding both in abstract and natural systems.
Hyperstructures apply to all kinds of universes: mathematical,
physical, chemical, biological, economic and social.  Furthermore, the
Transfer Principle makes it possible to connect them.  Detailed
applications will be the subject of future papers.  The main point in
this paper has been to illustrate the transfer of higher order binding
structure as given in a hyperstructure.  Ultimately one may also
consider bindings coming from hyperstructures of hyperstructures as in
the case of higher order links.

All the examples here and in Section \ref{sec:hyper} may be used to
put a hyperstructure on sets, spaces, structures and situations by the
methods described in Section \ref{sec:binding}.  This may be useful to
obtain actions like geometrical and physical fusion of objects in
various situations similar to the ``political/sociological'' metaphor.

After all we ``make things'' through a hyperstructure principle as in
modern engineering.  This may be so since it is the way nature works
through evolution, and after all we are ourselves products of such a
process.

\section{Multilevel state systems}
\label{sec:multi}
We have here discussed the organization of many-body systems or
general systems of collections of objects.  The systems may be finite,
infinite or even uncountable.  We have advocated Hyperstructures as
the guiding organizational principle.  In this section we will discuss
in more detail possible organization of the states of the system
through level connections.  We use the terminology in Section
\ref{sec:hyper}.

When we put a hyperstructure on a situation in order to obtain a
certain goal or action often a dynamics is required on $\calH$:
\begin{equation*}
  D \colon \calH \to \calH
\end{equation*}
which essentially changes the states.

In order to do this it is advantageous to be able to have as rich
structures as possible as states: sets, spaces (manifolds), algebras,
(higher) vector spaces, (higher) categories, etc.

For this reason and in order to cover as many interactions as possible
we introduce the following extension.

Instead of letting the states ($\Omega_i$) in $\calH$ take values in
$\sets$ we extend this to a family of prescribed \emph{hyperstructures
  of states}:
\begin{equation*}
  \calS = \{\calS_0,\calS_1,\ldots,\calS_n\}
\end{equation*}
where $\calS_i$ is a hyperstructure such that:
\begin{equation*}
  \begin{array}{c@{\quad}c@{\quad}l}
    \Omega_0 & \text{takes values in} & \calS_n\\
    \vdots&&\\
    \Omega_i & \text{takes values in} & \calS_{n - i}\\
    \vdots&&\\
    \Omega_n & \text{takes values in} & \calS_0.
  \end{array}
\end{equation*}
We want these level state structures to be connected in some way.
Therefore we require that $\calS$ is organized into a hyperstructure
itself with $\calS_0,\ldots,\calS_n$ as the levels --- or actually
sets of bonds of states, with $\calS_0$ being the top level (dual to
$\calH$ itself).  We furthermore assume that we have level connecting
assignments or boundary maps $\delta_i$ which we will for short write
as:
\begin{center}
  \begin{tikzpicture}
    \node (S0) at (0,0){$\calS_0$};
    \node (S1) at (1.5,0){$\calS_1$};
    \node (Sdots) at (3,0){$\cdots$};
    \node (Sn) at (4.5,0){$\calS_n$.};

    \begin{scope}[midway,above,font=\scriptsize]
      \draw[line join=round, decorate, decoration={
        zigzag,segment length=4,amplitude=.9,post=lineto,post
        length=2pt}, ->] (Sn) -- node{$\delta_n$} (Sdots);

      \draw[line join=round, decorate, decoration={
        zigzag,segment length=4,amplitude=.9,post=lineto,post
        length=2pt}, ->] (Sdots) -- node{$\delta_2$} (S1);

      \draw[line join=round, decorate, decoration={
        zigzag,segment length=4,amplitude=.9,post=lineto,post
        length=2pt}, ->] (S1) -- node{$\delta_1$} (S0);
    \end{scope}
  \end{tikzpicture}
\end{center}

Often the order of the hyperstructure will decrease moving from the
bottom to the top --- ``integrating away complexity''.  The
$\delta_i$'s may be of assignment (functional) or relational type.

This shows how to form hyperstructures of hyperstructures.  In such
cases one may actually use existing hyperstructures to form bonds and
states in new hyperstructures.

The assignment of a state to a bond (or collection of bonds) is a kind
of representation:
\begin{equation*}
  \Omega \colon \calH \rightsquigarrow \calS \quad \text{or} \quad
  \Omega_i \colon B_i \rightsquigarrow \calS_{n - i}
\end{equation*}
such that if
\begin{equation*}
  \Omega_i(b_i) = \omega_i \quad \text{and} \quad \partial_{i - 1} b_i
  = \{b_{i - 1}(\omega_{i - 1})\},
\end{equation*}
then
\begin{equation*}
  \delta_{n - i  + 1} \Omega_{i - 1} \partial_{i - 1} (b_i) = \Omega_i(b_i)
\end{equation*}
in simplified notation ($\{\text{-}\}$ meaning a family or subset of
objects).  This is a balancing equation of bonds and states at the
various levels.

Let us illustrate this by an example.

\begin{example}
  This is basically a version of Example d) in Section \ref{sec:hyper}
  and studied in \citeasnoun{HAM} in connection with genomic
  structure.  $\calH$ is given by sets
  \begin{equation*}
    G_0 < G_1 < G_2 < \cdots < G_n
  \end{equation*}
  meaning that there exist maps
  \begin{equation*}
    g_i \colon G_i \to \calP(G_{i + 1}).
  \end{equation*}
  To each $G_i$ we assign a state space $S_i$ --- possibly a manifold.

  Then the state hyperstructure reduces to the composition of (smooth)
  mappings:
  \begin{equation*}
    S_0 \leftarrow S_1 \leftarrow \cdots \leftarrow S_n.
  \end{equation*}
\end{example}

In order to influence global states from local actions it is a
reasonable and general procedure to put a hyperstructure on the
systems with a multilevel state structure given by another
hyperstructure $\calS$ with level relations as described:
\begin{center}
  \begin{tikzpicture}
    \node (S0) at (0,0){$\calS_0$};
    \node (S1) at (1.5,0){$\calS_1$};
    \node (Sdots) at (3,0){$\cdots$};
    \node (Sn) at (4.5,0){$\calS_n$.};

    \draw[line join=round,decorate, decoration={
      zigzag,segment length=4,amplitude=.9,post=lineto,post
      length=2pt}, ->] (Sn) -- (Sdots);

    \draw[line join=round,decorate, decoration={
      zigzag,segment length=4,amplitude=.9,post=lineto,post
      length=2pt}, ->] (Sdots) -- (S1);

    \draw[line join=round,decorate, decoration={
      zigzag,segment length=4,amplitude=.9,post=lineto,post
      length=2pt}, ->] (S1) -- (S0);
  \end{tikzpicture}
\end{center}
The idea is then to act on $\calS_n$ by introducing a suitable
dynamics and let the actions propagate through the hyperstructure to
the global level.  This is similar to social systems and may be called
the ``Democratic Method of Action''.  In other situations one may want
to move from high to low level states.

It is especially useful if the level relations are functional
assignments:
\begin{equation*}
  \calS_0 \leftarrow \calS_1 \leftarrow \cdots \leftarrow \calS_n.
\end{equation*}
For example if the $\calS_i$'s are categories of some kind, the arrows
would represent functors and if the $\calS_i$'s are spaces, the arrows
represent mappings.

Let us specify the mappings
\begin{equation*}
  \partial_{i - 1} \colon X_i \to \calP(X_{i - 1})
\end{equation*}
by
\begin{equation*}
  \partial_{i - 1} b_i = \{b_{i - 1}(\omega_{i - 1})\},
\end{equation*}
where $\{-\}$ just means a family or subsets of objects, and requiring
assignments (mappings):
\begin{center}
  \begin{tikzpicture}
    \node (prod) at (0,0){$\prod \Omega_{i - 1}(\{b_{i - 1}\})$};
    \node (omega) at (3,0){$\Omega_i(b_i)$};
    \node (elementinprod) at (0,-0.5){\rotatebox{-90}{$\in$}};
    \node (elementinomega) at (3,-0.5){\rotatebox{-90}{$\in$}};
    \node (Sn-i+1) at (0,-1){$\calS_{n - i + 1}$};
    \node (Sn-i) at (3,-1){$\calS_{n - i}$};

    \draw[line join=round,decorate, decoration={
      zigzag,segment length=4,amplitude=.9,post=lineto,post
      length=2pt}, ->] (prod) -- node[above, midway,
        font=\scriptsize]{$\hat{\delta}_{n - i + 1}$} (omega);

    \draw[line join=round,decorate, decoration={
      zigzag,segment length=4,amplitude=.9,post=lineto,post
      length=2pt}, ->] (Sn-i+1) -- node[above, midway,
        font=\scriptsize]{$\delta_{n - i + 1}$} (Sn-i);
  \end{tikzpicture}
\end{center}
The $\hat{\delta}$'s are then state level connectors.  If the
$\calS_i$'s have a tensor product we will often require
\begin{equation*}
  \bigotimes\{\Omega_{i - 1}(b_{i - 1})\} \xrightarrow{\hat{\delta}_{n
      - i + 1}} \Omega_i(b_i)
\end{equation*}
%
More schematically this may be described as:
\begin{equation*}
  \{b(i_n)\} \xrightarrow{\partial} \{b(i_{n - 1},i_n)\}
  \xrightarrow{\partial} \cdots \xrightarrow{\partial}
  \{b(i_0,\ldots,i_n)\}
\end{equation*}
for bonds, and for states:
\begin{equation*}
  \Omega_n(\{b(i_n)\}) \xleftarrow{\hat{\delta}} \Omega_{n -
    1}(\{b(i_{n - 1},i_n)\}) \xleftarrow{\hat{\delta}} \cdots
  \xleftarrow{\hat{\delta}} \Omega_0(\{b(i_0,\ldots,i_n)\}).
\end{equation*}
This is an extension of the framework for extended Topological Quantum
Field Theories (see \citeasnoun{Lawrence}) where one considers
manifolds with general boundaries and corners
\begin{equation*}
  M = \{M(\omega)\} \quad \text{like in Section \ref{sec:app} i)}
\end{equation*}
bound by generalized cobordisms.

At the state level one assigns to these manifolds higher order
algebraic structures such as a version $k$-vector spaces ($k-\calV$)
or $k$-algebras as follows \cite{Lawrence}:
\begin{equation*}
  \begin{array}{r@{\ }c@{\ }l}
    M(\omega) & \longrightarrow & Z(\omega) \in k-\calV\\
    (\codim = k) && (\order = k)
  \end{array}
\end{equation*}
where levels of states are connected via geometrically induced
pairings:
\begin{equation*}
  k-\calV \otimes k-\calV \to (k - 1)-\calV
\end{equation*}
for all $k$.  $0-\calV$ being the scalars $\C$ in the case of complex
vector spaces.

In this sense we can control and regulate the global state from the
lowest level which is clearly a desirable thing in systems of all
kinds.  This is useful in the following situation.  Let $A$ be a
desired action or task which may be ``un-managable'' in a given system
or context.  Furthermore, let $M = \{m_\alpha\}$ be a collection of
``managable'' actions in the system.

Put $X_0 = M$ and design a hyperstructure $\calH$ (like a society or
factory) where then $A$ appears as $X_n = B_n$ --- the top level bond
of the hyperstructure $\calH$.  Hence $\calH$ will act as a propagator
from \{managable actions\} to \{desired actions\} by dynamically
regulating the states of $\calH$ as described.  This procedure applies
to general systems and in the sense of \citeasnoun{Waddington} one may
say that hyperstructures are ``Tools for Thought'' and creation of
novelty.

In general systems one may often want to change the global state in a
desired way and it may be difficult since it would require large
resources (``high energy'').  But via a hyperstructure it may be
possible introducing managable local actions and change of states
which may require small resources (``low energy''), in other words
small actions are being organized into large actions via $\calH$.
This is similar to changing the action of a society or organization by
influencing individuals.  If one wants to join two opposing societies
(or nations) into one, it may take less resources to act on
individuals to obtain the global effect.  The photosynthesis works
along the same lines --- collecting, organizing and amplifying
energy.

It seems like an interesting idea to suggest the use of
hyperstructures in order to facilitate fusion of various types of
systems, for example ``particles'' in biology, chemistry and physics.
Even nuclear fusion may profit from this perspective.  The
hyperstructure in question may be introduced on the system itself or
surrounding space and forces.

\section{Organizing principles}
\label{sec:org}
Binding structures and hyperstructures as we have described them are
basically organizing principles of collections of objects.  They apply
to all kinds of collections and general systems, and as organizing
principles they may be particularly interesting in physical matter
(condensed matter) of atoms, moleclues, etc.

R.\ Laughlin has advocated the importance of organizing principles in
condensed matter physics in understanding for example
superconductivity, the quantum Hall effect, phonons, topological
insulators etc.  See \citeasnoun{LP}.  He suggests that the very
precise measurements made of important physical constants in these
situations come from underlying organizing principles of matter.

Our binding- and hyperstructures are organizing principles that when
introduced to physical matter should lead to new emergent properties.
It is like in our Example \ref{ex:links}.  If we are given a random
tangle of links, a new and non-trivial geometric order emerges when we
put a higher order Brunnian structure on the collection of links.\\

In one way it is analogous to logical systems, where organized
statements are more likely to be decidable than random ones.
Similarly in biology: language, memory, spatial recognition, etc.\
are related to similar organizing principles.

Entangled states are studied in quantum mechanics and higher order
versions are suggested \citeasnoun{Z}.  Greenberger--Horne--Zeilinger
(GHZ) states are analogous to Borromean and Brunnian rings, see
\citeasnoun{A}.  From our previous discussion we are naturally led to
suggest higher order entangled states organized by a hyperstructure
$\calH$ using the transfer principle.  Could \ such \ a \ process \
lead \ to \ collections \ of \ particles \ forming \
global/macroscopic quantum states?  Could also an $\calH$-structure
act as a kind of (geometric) protectorate of a desired quantum state
from thermodynamic disturbances (in high temperature
superconducitivity for example)?

The binding principle may be applied in two ways --- in particular in
condensed matter physics.\\

\begin{itemize}
\item[I.] Putting a binding- or hyperstructure on a collection
  of objects.  Then collections of bound structures will appear, and
  they may have interesting emergent properties.  For example with
  respect to precise measurements of involved constants of nature.\\
\item[II.] Putting a binding- or hyperstructure on the ambient space
  (space-time) of a collection for example using various fields etc.
  This will introduce a structure on the collection and may result in
  bindings and fusion of particles and objects, or splitting
  (fission), stabilizing them into new patterns with new emergent
  properties.\\
\end{itemize}

In other words space, fields and reactors may all be organized by
binding principles.  We suggest that putting a binding- or
hyperstructure on a collection, situation or system is a very
fundamental and useful organizing principle.

\subsection*{Acknowledgements}
I would like to thank A.~Stacey and M.~Thaule for help with graphics
and technical issues.  All the figures appear with their kind
permission.

\subsection*{Related links of interest:}
\begin{itemize}
  \item \href{http://www.newscientist.com/article/mg20927942.300-make-way-for-mathematical-matter.html}{New
    Scientist}\footnote{\href{http://www.newscientist.com/article/mg20927942.300-make-way-for-mathematical-matter.html}{\texttt{http://www.newscientist.com/article/mg20927942.300-make-way-for-
          mathematical-matter.htm}}}
  \item
    \href{http://www.technologyreview.com/view/422055/topologist-predicts-new-form-of-matter/}{MIT
      Tech
      Review}\footnote{\href{http://www.technologyreview.com/view/422055/topologist-predicts-new-form-of-matter/}{\texttt{http://www.technologyreview.com/view/422055/topologist-predicts-
          new-form-of-matter/}}}
\end{itemize}

\end{document}